\documentclass{article}

\usepackage[giveninits=true,maxbibnames=9,maxcitenames=2,backend=biber,url=false,doi=false,isbn=false]{biblatex}

    \usepackage[final, nonatbib]{neurips_2025}

\usepackage[utf8]{inputenc} %
\usepackage[T1]{fontenc}    %
\usepackage{hyperref}       %
\usepackage{url}            %
\usepackage{booktabs}       %
\usepackage{amsfonts}       %
\usepackage{nicefrac}       %
\usepackage{microtype}      %
\usepackage{xcolor}         %

\usepackage{amsmath,amsthm,epsfig,amssymb,amsbsy}
\usepackage{enumitem}
\usepackage{graphicx}

\DeclareMathOperator*\minimize{minimize}

\DeclareMathOperator\stt{subject\ to}

\DeclareMathOperator{\dom}{dom}
\DeclareMathOperator{\diag}{diag}

\DeclareMathOperator{\intr}{int}

\DeclareMathOperator{\sign}{sign}

\DeclareMathOperator{\barsgn}{\overline{sgn}}
\DeclareMathOperator{\arcsinh}{arcsinh}
\DeclareMathOperator{\trace}{trace}
\DeclareMathOperator{\laspan}{span}
\DeclareMathOperator{\vecop}{vec}

\newcommand{\bR}{\mathbb{R}}

\newcommand{\bB}{\mathbb{B}}

\newcommand{\bP}{\mathbb{P}}

\newcommand{\bS}{\mathbb{S}}
\newcommand{\exR}{\overline{\mathbb{R}}}

\newcommand{\cA}{\mathcal{A}}
\newcommand{\cF}{\mathcal{F}}
\newcommand{\cC}{\mathcal{C}}

\newcommand{\cG}{\mathcal{G}}

\newcommand{\cX}{\mathcal{X}}
\newcommand{\cO}{\mathcal{O}}
\newcommand{\cZ}{\mathcal{Z}}
\newcommand{\cT}{\mathcal{T}}
\newcommand{\cH}{\mathcal{H}}

\newcommand\R{\mathbb{R}}
\newcommand\N{\mathbb{N}}
\newcommand\C{\mathcal{C}}

\newcommand\seq[1]{(#1)_{k\in\N}}

\newcommand{\Cpp}{C\kern-0.04em+\kern-0.01em+}

\DeclareMathOperator{\proj}{\Pi}

\usepackage[]{cleveref}[0.19]

\newtheorem{theorem}{Theorem}[section]
\makeatletter
\newcommand{\settheoremtag}[1]{%
  \let\oldthetheorem\thetheorem%
  \renewcommand{\thetheorem}{#1}%
  \g@addto@macro\endtheorem{%
    \addtocounter{theorem}{-1}%
    \global\let\thetheorem\oldthetheorem}%
  }
\makeatother

\newtheorem{lemma}[theorem]{Lemma}
\newlist{lemenum}{enumerate}{1} %
\setlist[lemenum]{label=(\roman*), ref=\theproposition(\roman*), font=\rm}
\crefalias{lemenumi}{Lemma} 

\newtheorem{proposition}[theorem]{Proposition}
\newlist{propenum}{enumerate}{1} %
\setlist[propenum]{label=(\roman*), ref=\theproposition(\roman*), font=\rm}
\crefalias{propenumi}{Proposition} 

\newlist{theoremenum}{enumerate}{1} %
\setlist[theoremenum]{label=(\roman*), ref=\thetheorem(\roman*), font=\rm}
\crefalias{theoremenumi}{Theorem} 

\newtheorem{definition}[theorem]{Definition}

\crefname{theorem}{Theorem}{Theorems}
\crefname{Theorem}{Theorem}{Theorems}

\newtheorem{assumption}[theorem]{Assumption}
\crefname{assumption}{Assumption}{Assumptions}
\Crefname{assumption}{Assumption}{Assumptions}

\newtheorem{example}[theorem]{Example}

\usepackage{algorithm}
\usepackage{algpseudocode}

\newcommand\algorithmicinitialize{\textsc{Initialize:}}
\algnewcommand\Initialize{\item[\algorithmicinitialize]}%

\title{Nonlinear preconditioning for optimization under generalized smoothness}

\title{A nonlinear preconditioning perspective on generalized smoothness in optimization}

\title{Escaping saddle points without Lipschitz smoothness: the power of nonlinear preconditioning}

\author{%
Alexander Bodard\\
KU Leuven\\
\texttt{alexander.bodard@kuleuven.be}\\
\And
Panagiotis Patrinos\\
KU Leuven\\
\texttt{panos.patrinos@esat.kuleuven.be}\\
}

\begin{document}

\maketitle

\begin{abstract}
  We study generalized smoothness in nonconvex optimization, focusing on \((L_0, L_1)\)-smoothness and anisotropic smoothness. The former was empirically derived from practical neural network training examples, while the latter arises naturally in the analysis of nonlinearly preconditioned gradient methods. We introduce a new sufficient condition that encompasses both notions, reveals their close connection, and holds in key applications such as phase retrieval and matrix factorization. Leveraging tools from dynamical systems theory, we then show that nonlinear preconditioning -- including gradient clipping -- preserves the saddle point avoidance property of classical gradient descent. Crucially, the assumptions required for this analysis are actually satisfied in these applications, unlike in classical results that rely on restrictive Lipschitz smoothness conditions. We further analyze a perturbed variant that efficiently attains second-order stationarity with only logarithmic dependence on dimension, matching similar guarantees of classical gradient methods. 

\end{abstract}

\section{Introduction} \label{sec:intro}

We consider the unconstrained optimization problem
\begin{equation}
    \minimize_{x \in \mathbb{R}^n} f(x),
\end{equation}
where \( f: \mathbb{R}^n \to \mathbb{R} \) is a twice continuously differentiable nonconvex function. This work studies the \emph{nonlinearly preconditioned gradient method}, with iterates described by
\begin{equation} \label{eq:pgd} \tag{P-GD}
    x^{k+1} = T_{\gamma, \lambda}(x^k) := x^k - \gamma \nabla \phi^*(\lambda \nabla f(x^k)),
\end{equation}
where \(\phi : \R^n \to \R \cup \{\infty\}\) is referred to as the \emph{reference function}, and \(\phi^*\), its convex conjugate, is called the \emph{dual reference function}.

\begin{figure}[t]
    \centering
        \includegraphics[width=0.46\linewidth]{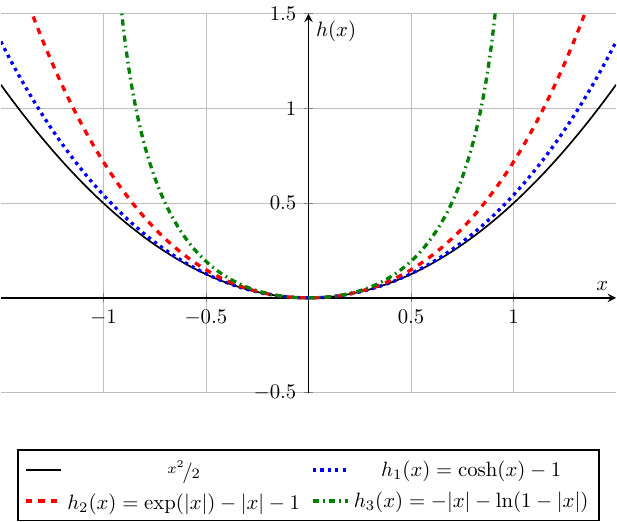}
    \hfill
        \includegraphics[width=0.46\linewidth]{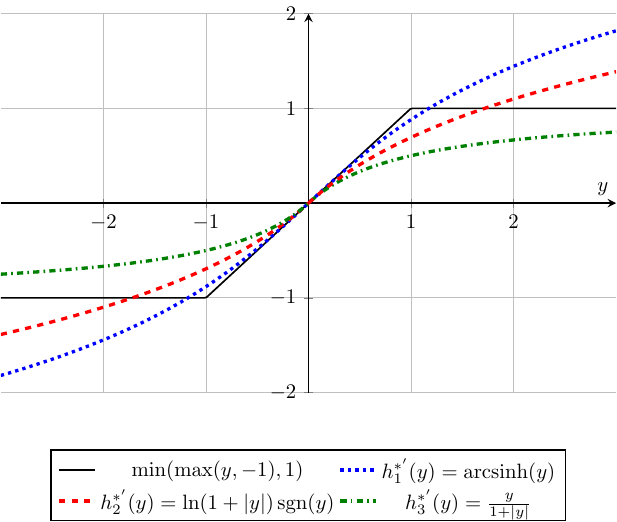}
    \caption{Comparison of (a) kernel functions and (b) their corresponding nonlinear preconditioners.}
    \label{fig:kernels}
\end{figure}

Nonlinear preconditioning provides a flexible framework for constructing and analyzing gradient-based optimization algorithms \cite{oikonomidis_nonlinearly_2025,laude_anisotropic_2025,maddison_dual_2021}. For instance, when \( \phi(x) = \frac{1}{2} \Vert x \Vert^2 \), the update \eqref{eq:pgd} reduces to classical gradient descent. More broadly, we focus on \emph{isotropic} reference functions of the form \( \phi(x) = h(\Vert x \Vert) \) for some scalar kernel function \( h : \mathbb{R} \to \mathbb{R}_+ \cup \{ \infty \} \), though our results extend in part to more general settings, including \emph{separable} reference functions \( \phi(x) = \sum_{i=1}^n h(x_i) \).
Some kernel functions of interest include:
\begin{equation} \label{eq:h-examples}
    h_1(x) = \cosh(x) - 1, \quad h_2(x) = \exp(\vert x \vert) - \vert x \vert - 1, \quad h_3(x) = - \vert x \vert - \ln(1 - \vert x \vert),
\end{equation}
each of which upper bounds the quadratic function \( \nicefrac{x^2}{2} \), as visualized in Fig.\;\ref{fig:kernels}. These choices induce preconditioners that closely resemble common \emph{gradient clipping} heuristics, as shown in Fig.\;\ref{fig:kernels}. 

The effectiveness of gradient clipping has been justified using the concept of \((L_0, L_1)\)-smoothness, which is empirically motivated by practical neural network training scenarios \cite{zhang_why_2020}. However, it remains unclear under what precise conditions this smoothness assumption holds in real-world applications.
On the other hand, the preconditioned gradient method is naturally analyzed under \emph{anisotropic smoothness} \cite{oikonomidis_nonlinearly_2025},
another generalization of the classical Lipschitz smoothness condition. 
Rather than imposing a global \emph{quadratic} upper bound, anisotropic smoothness permits more flexible upper bounds defined in terms of the reference function \( \phi \). This makes the preconditioned gradient method particularly attractive in settings where the standard Lipschitz condition is too restrictive.
This leads us to our first central question:

\medskip
\noindent
\emph{Can we formally establish anisotropic smoothness and \((L_0, L_1)\)-smoothness of practical problems where traditional assumptions fail?}

\medskip

Our second line of inquiry focuses on the behavior of the preconditioned gradient method when applied to nonconvex objectives. Classical gradient descent is known to \emph{avoid strict saddle points} under the assumption of (global) Lipschitz smoothness \cite{lee_first-order_2019}, a phenomenon which helps explain its strong empirical performance in nonconvex settings.
However, for many practical applications Lipschitz smoothness holds only \emph{locally} or on compact sets around a minimizer, meaning that this assumption is not truly satisfied.
This raises the following question:

\medskip
\noindent
\emph{Does nonlinear preconditioning preserve the saddle point avoidance properties of gradient descent under a possibly less stringent smoothness assumption?}

\medskip

Our results reveal novel connections between different generalizations of smoothness and provide strong theoretical support for nonlinear preconditioning, particularly in nonconvex settings where the classical Lipschitz smoothness assumption may fail.

\paragraph{Contributions}

Our contributions can be summarized as follows.

\begin{itemize}
    \item We investigate the classes of problems for which \((L_0, L_1)\)-smoothness and anisotropic smoothness -- two generalizations of the classical Lipschitz smoothness condition -- are applicable. To this end, we propose a novel sufficient condition (\cref{assump:novel-condition}) that guarantees both anisotropic and \((L_0, L_1)\)-smoothness, thereby revealing a structural link between these two frameworks. We further demonstrate in \cref{sec:applications} that this condition holds for several prominent nonconvex problems, including phase retrieval, low-rank matrix factorization, and Burer-Monteiro factorizations of MaxCut-type problems.
    \item We establish that nonlinear preconditioning preserves the saddle point avoidance behavior of gradient descent, and moreover extends results from the classical Lipschitz smoothness framework to the broader setting of anisotropic smoothness. 
    Specifically, we prove asymptotic avoidance of strict saddle points by leveraging the stable-center manifold theorem. 
    By invoking a recent nonsmooth generalization of this theorem, this analysis is then further extended to accommodate hard gradient clipping. 
    Finally, we present a complexity analysis for a perturbed variant of the preconditioned gradient method, showing that it converges to a second-order stationary point with only logarithmic dependence on the problem dimension.

\end{itemize}

\paragraph{Notation}

Let \(\bS^{n \times n}\) be the set of symmetric \(n \times n\) matrices.
We denote the standard Euclidean inner product on \(\R^n\) by \(\langle \cdot, \cdot, \rangle\), and the corresponding norm by \(\Vert \cdot \Vert\).
For \(X, Y \in \R^{m \times n}\), \(\langle X, Y \rangle = \trace(X^\top Y)\) is the standard inner product on \(\R^{m \times n}\) and \(\Vert \cdot \Vert\) denotes the spectral norm.
The class of \(k\) times continuously differentiable functions on an open set \(O \subseteq \R^n\) is denoted by \(\cC^k(O)\). 
We write \(\barsgn(x) = \nicefrac{x}{\Vert x \Vert}\) for \(x \in \R^n \setminus \{0\}\) and \(0\) otherwise.
A function \(f \in \cC^2(\R^n)\) is \(L\)-Lipschitz smooth if for all \(x, y \in \R^n\) it holds that \(\Vert \nabla f(x) - \nabla f(y) \Vert \leq L \Vert x - y \Vert\), with \(L \geq 0\), and \((L_0, L_1)\)-smooth if \(\Vert \nabla^2 f(x) \Vert \leq L_0 + L_1 \Vert \nabla f(x) \Vert \) for all \(x \in \R^n\) with \(L_0, L_1 \geq 0\). 
Otherwise, we follow \cite{rockafellar_variational_1998}.

\subsection{Related work}

\paragraph{Generalized smoothness} 
Gradient descent is traditionally analyzed under the assumption of Lipschitz smoothness~\cite{nesterov_lectures_2018}, although many applications violate this condition. Bregman relative smoothness is a popular extension which allows the Hessian to grow unbounded, see e.g.\,\cite{lu_relatively_2018} which assumes a certain polynomial growth.
More recently, the \((L_0, L_1)\)-smoothness condition was proposed by Zhang et al.~\cite{zhang_why_2020}, based on empirical observations in LSTMs, and used to analyze clipped gradient descent and a momentum variant~\cite{zhang_improved_2020}. The framework has since been applied to stochastic normalized gradient descent~\cite{zhao_convergence_2021} and generalized SignSGD~\cite{crawshaw_robustness_2022}. Notably, Crawshaw et al.~\cite{crawshaw_robustness_2022} provided empirical evidence that \((L_0, L_1)\)-smoothness holds for Transformers~\cite{vaswani_attention_2017}, albeit with layer-wise variation in constants.
Further generalizations include \(\alpha\)-symmetric smoothness~\cite{chen_generalized-smooth_2023} and \(\ell\)-smoothness~\cite{li_convex_2023}, and the latter was used to analyze the convergence of Adam~\cite{li_convergence_2023}. 
Despite empirical support for these conditions in key applications, theoretical guarantees remain limited.

\paragraph{Nonlinear preconditioning} 
The preconditioned gradient method with updates given by~\eqref{eq:pgd} was introduced in the convex setting by Maddison et al.~\cite{maddison_dual_2021}. Then, Laude et al.~\cite{laude_dualities_2023,laude_anisotropic_2025} studied \(L\)-anisotropic smoothness and, under this condition, showed convergence of \eqref{eq:pgd} for nonconvex problems. 
The method was later extended to measure spaces~\cite{bonet_mirror_2024}. Oikonomidis et al.~\cite{oikonomidis_nonlinearly_2025} proposed the \((L, \bar{L})\)-anisotropic smoothness condition, connected it to \((L_0, L_1)\)-smoothness, and analyzed convergence of~\eqref{eq:pgd} in both convex and nonconvex settings.
We also highlight the works \cite{leger_gradient_2023,oikonomidis_forward-backward_2025} that study the concept of \(\Phi\)-convexity, which is closely related to anisotropic smoothness.

\paragraph{Saddle point avoidance} 
To explain the success of gradient descent on \textit{nonconvex} problems, much work has focused on its (strict) saddle point avoidance properties~\cite{lee_gradient_2016, lee_first-order_2019}. It was shown that gradient descent may take exponential time to escape saddle points, even with random initialization~\cite{du_gradient_2017}. 
The works \cite{levy_power_2016,murray_revisiting_2019} showed that noise-injected normalized gradient descent escapes them more efficiently.
Jin et al.~\cite{jin_how_2017,jin_nonconvex_2021} demonstrated that perturbed gradient descent escapes saddle points in time polylogarithmic in the problem dimension. Recently, Cao et al.~\cite{cao_efficiently_2025} studied saddle point avoidance under a second-order self-bounding regularity condition rather than under classical Lipschitz smoothness.

\section{Anisotropic smoothness} \label{sec:anisotropic-smoothness}
\subsection{Definition and basic properties}

This section introduces \((L, \bar L)\)-anisotropic smoothness as proposed by \cite{oikonomidis_nonlinearly_2025}, and summarizes some known properties which are relevant to this work.
The following assumption, which guarantees in particular that \(\phi^* \in \cC^1(\R^n)\) and \(\phi \geq 0\), is considered valid throughout.
\begin{assumption} \label{assump:basic-reference}
    The function \(\phi : \R^n \to \exR\) is proper, lsc, strongly convex and even with \(\phi(0) = 0 \).
\end{assumption}
We usually also assume the following condition, which ensures in particular that \(\phi^* \in \cC^2(\R^n)\).
\begin{assumption} \label{assumption:phi-star-c2}
    \( \intr \dom \phi \neq \emptyset \); \( \phi \in \cC^2(\intr \dom \phi) \), and for any sequence \(\{x^k\}_{k \in \N}\) that converges to some boundary point of \(\intr \dom \phi\), it follows that \( \Vert \nabla \phi(x^k) \Vert \to \infty \).
\end{assumption}
We follow the definition of anisotropic smoothness by \cite{oikonomidis_nonlinearly_2025}, which reduces to \cite[Definition 3.1]{laude_anisotropic_2025} with reference function \(\bar L \phi\) if \(\dom \phi = \R^n\). 
For a geometric intuition, we refer to \cite{leger_gradient_2023,oikonomidis_forward-backward_2025,oikonomidis_nonlinearly_2025}.
\begin{definition}[\((L, \bar L)\)-anisotropic smoothness {\cite{oikonomidis_nonlinearly_2025}}] \label{def:anisotropic-smoothness}
    A function \( f : \R^n \to \R \) is \((L, \bar L)\)-anisotropically smooth relative to a reference function \( \phi \) with constants \(L, \bar L > 0\) if
    \begin{equation*}
        f(x) \leq f(\bar x) + \bar L L^{-1}\phi(L(x - \bar y)) - \bar L L^{-1} \phi(L(\bar x - \bar y))
    \end{equation*}
    for all \(x, \bar x \in \R^n\), where \( \bar y = T_{L^{-1}, \bar L^{-1}} (\bar x) = \bar x - L^{-1} \nabla \phi^*(\bar L^{-1} \nabla f(\bar x)) \).
\end{definition}
The following proposition provides a sufficient condition for anisotropic smoothness.
We consider the case \(\phi^* \in \cC^2\) for simplicity of exposition, but note that a variant for \(\phi^* \notin \cC^2\) can also be formulated.
\begin{proposition}[Second-order characterization of \((L, \bar L)\)-anisotropic smoothness] \label{prop:second-order-characterization}
    Suppose that \cref{assumption:phi-star-c2} holds, and let \(f \in \cC^2\) be such that for all \(x \in \R^n\)
    \begin{equation} \label{eq:second-order-characterization}
        \lambda_{\max} (\nabla^2 \phi^*(\bar L^{-1} \nabla f(x)) \nabla^2 f(x)) \leq L \bar L,
    \end{equation}
    and \(\lim_{\Vert x \Vert \to \infty} \Vert T_{L^{-1}, \bar L^{-1}}(x) \Vert = \infty\).
    Moreover, assume that either \( \dom \phi \) is bounded or that \( \dom \phi = \R^n\), and that for all \( x \in \R^n \) we have
    \(
        f(x) \leq \bar L r^{-1} \phi(r x) - \beta
    \)
    for some \(r \in (0, L)\), \( b \in \R\).
    Then, \(f\) is \((\delta L, \bar L)\)-anisotropically smooth relative to \(\phi\) for any \(\delta > 1\).
\end{proposition}
We say that \(f\) satisfies the \emph{second-order characterization} of anisotropic smoothness if \eqref{eq:second-order-characterization} holds.
Note that the growth condition on \(f\) is not restrictive when \(\phi = \dom \R^n\), and that the coercivity assumption on the iteration map \(T_{L^{-1}, \bar L^{-1}}\) is very mild; we refer the reader to the arguments in \cite{oikonomidis_nonlinearly_2025}.
Finally, we connect anisotropic smoothness to some popular smoothness notions.
\begin{example}[Lipschitz-smoothness {\cite[Proposition 2.3]{oikonomidis_nonlinearly_2025}}]
    Suppose that \(f \in \cC^2\) is \(L_f\)-Lipschitz smooth.
    Denote by \(\mu > 0\) the parameter of strong convexity of a reference function \(\phi\).
    Then \(f\) is \((\nicefrac{L_f}{\mu}, 1)\)-anisotropically smooth relative to \(\phi\). 
\end{example}

\begin{example}[\((L_0, L_1)\)-smoothness] \label{ex:l0-l1-smoothness}
    Let \(f \in \cC^2\) be \((L_0, L_1)\)-smooth, let \(L = L_1, \bar L = \nicefrac{L_0}{L_1}\), and let \(\phi(x) = -\Vert x \Vert - \ln(1 - \Vert x \Vert)\).
    Then \(f\) satisfies the second-order characterization of \((L, \bar L)\)-anisotropic smoothness relative to \(\phi\) \cite[Proposition 2.6 \& Corollary 2.7]{oikonomidis_nonlinearly_2025}.
\end{example}

\subsection{A novel sufficient condition for generalized smoothness}

Although it is well-known that univariate polynomials are \((L_0, L_1)\)-smooth \cite[Lemma 2]{zhang_why_2020}, this is not necessarily the case for multivariate polynomials, as illustrated by the following example.
\begin{example} \label{ex:multivariate}
    Consider the polynomial \(f(x, y) = \frac{1}{4} x^4 + \frac{1}{4}y^4 -\frac{1}{2} x^2 y^2\) with gradient and Hessian
    \begin{equation*}
        \nabla f(x, y) = \begin{pmatrix}
            x^3 - xy^2\\
            y^3 - x^2 y
        \end{pmatrix}, \qquad \nabla^2 f(x, y) = \begin{pmatrix}
            3x^2 - y^2 & -2xy\\
            -2xy & -x^2 + 3 y^2
        \end{pmatrix}.
    \end{equation*}
    Remark that \(\nabla f(x, -x) = 0\) and \(\nabla^2 f(x, -x) = x^2 \left( \begin{smallmatrix}
        2 & -2\\
        - 2 & 2
    \end{smallmatrix} \right)\).
    Clearly, \(f\) cannot be \((L_0, L_1)\)-smooth since \(\Vert \nabla^2 f(x, -x) \Vert_F = 4 \Vert x \Vert^2\) grows unbounded, while \(\Vert \nabla f(x, -x) \Vert = 0\) for all \(x \in \R\). 
\end{example}
For multivariate polynomials there may exist a path of \(\Vert x \Vert \to \infty\) along which the gradient norm grows \emph{slower} than the Hessian norm, in which case \((L_0, L_1)\)-smoothness cannot hold.
More examples are included in \cref{sec:multivariate-polynomials-fail-extended}.
Based on this insight, we propose the following  novel condition.  
\begin{assumption} \label{assump:novel-condition}
    There exists an \(R \in \N\) such that for all \(x \in \R^n\)
    \begin{equation*}
        \Vert \nabla^2 f(x) \Vert_F \leq p_{R}(\Vert x \Vert), \quad \text{and} \quad \Vert \nabla f(x) \Vert \geq q_{R+1}(\Vert x \Vert).
    \end{equation*}
    Here \(p_R(\alpha) = \sum_{i = 0}^R a_i \alpha^i\) and \(q_{R+1}(\alpha) = \sum_{i = 0}^{R+1} b_i \alpha^i\) are polynomials of degree \(R\) and \(R+1\), respectively, and in particular we assume that \(b_{R+1} > 0\).
\end{assumption}
Note that \cite{lu_relatively_2018} constructs a Bregman distance inducing kernel function under a similar polynomial upper bound to the Hessian norm.
\Cref{sec:univariate-polynomials-verification} verifies \cref{assump:novel-condition} for univariate polynomials.
The following result states that \cref{assump:novel-condition} is a sufficient condition for \((L_0, L_1)\)-smoothness.\footnote{
        In fact, \cref{thm:l0-l1-smoothness-polynomials} still holds if \cref{assump:novel-condition} is relaxed to \(\Vert \nabla f(x) \Vert \geq q_{R}(\Vert x \Vert)\).
}
\begin{theorem} \label{thm:l0-l1-smoothness-polynomials}
    Suppose that \cref{assump:novel-condition} holds for \(f \in \cC^2\).
    Then, for any \( L_1 > 0\) there exists an \( L_0 > 0\) such that \(f\) is \((L_0, L_1)\)-smooth.
\end{theorem}
Under mild conditions on the kernel function \(h\), which \cref{sec:assumptions} shows hold for all examples in \eqref{eq:h-examples}, \cref{assump:novel-condition} also implies the second-order characterization of anisotropic smoothness.
In fact, it implies the stronger condition that \(\Vert \nabla^2 \phi^*(\bar L^{-1} \nabla f(x)) \nabla^2 f(x) \Vert\) is uniformly bounded.
\begin{assumption} \label{assumption:isotropic}
    The reference function \(\phi\) is isotropic, i.e., \(\phi(x) = h(\Vert x \Vert)\), and such that \normalfont{(i)} \( \nicefrac{{h^*}'(y)}{y} \) is a decreasing function on \( \R_+\), \normalfont{(ii)} \(\lim_{y \to +\infty} y {h^*}''(y) = C_2\), for some \(C_2 \in \R_+\), and \normalfont{(iii)}
    \begin{equation*}
        \begin{aligned}
            &\lim_{y \to +\infty} \frac{{h^*}'(s_d(y))}{y} = 0, \qquad \text{for any polynomial } s_{d}(\alpha) = \sum_{i = 0}^d u_i \alpha^i \text{ of degree } d.
        \end{aligned}
    \end{equation*}
\end{assumption}
\begin{theorem} \label{thm:anisotropic-smoothness-polynomials}
    Suppose that \(f\) satisfies \cref{assump:novel-condition}.
    If \(\phi\) satisfies \cref{assumption:phi-star-c2} and \cref{assumption:isotropic}, then for any \( \bar L > 0\) there exists an \(L > 0\) such that \(f\) satisfies the second-order characterization of \((L, \bar L)\)-anisotropic smoothness relative to \(\phi\). 
\end{theorem}

\subsection{Applications} \label{sec:applications}

We now establish for a number of key applications that \cref{assump:novel-condition} holds, thus proving that the objective is \((L_0, L_1)\)-smooth and satisfies the second-order characterization of \((L, \bar L)\)-anisotropic smoothness. 
Remark that for all of these, the classical Lipschitz smoothness assumption is violated.

\subsubsection{Phase retrieval}

Consider the real-valued phase retrieval problem with objective and gradient
\begin{equation} \label{eq:phase-retrieval-objective-gradient}
    f(x) = \frac{1}{4} \sum_{i = 1}^m \left( y_i^2 - (a_i^\top x)^2 \right)^2, \quad \nabla f(x) = -\sum_{i = 1}^m \left( y_i^2 - (a_i^\top x)^2 \right) a_i a_i^\top x. %
\end{equation}
Here, \(a_i \in \R^n\) and \(y_i \in \R\) for \(i \in \N_{[1, m]}\) are the measurement vectors and the corresponding measurements, respectively.
A relaxed smoothness condition for the phase retrieval problem has been explored in \cite{bolte2018first} based on Bregman distances.
The following theorem establishes that whenever the measurement vectors span \(\R^n\), the objective \(f\) also satisfies our \cref{assump:novel-condition}.
Note that the measurement vectors can only span \(\R^n\) if \(m \geq n\).
Moreover, the assumption on spanning \(\R^n\) is mild compared to well-studied conditions that guarantee signal recovery in the phase retrieval problem.
These conditions either require randomly sampled measurement vectors with \(m\) on the order of \(n \log n\) \cite{candes_phaselift_2013}, or the so-called complement property \cite{bandeira_saving_2014}.
The former ensures the spanning property with high probability, while the latter guarantees it deterministically.
\begin{theorem} \label{th:phase-retrieval-smoothness}
    Consider the phase retrieval problem with objective \eqref{eq:phase-retrieval-objective-gradient} and suppose that the vectors \(\{a_i\}_{i = 1}^m\) span \(\R^n\).
    \begin{theoremenum}
        \item \label{th:phase-retrieval-smoothness-1} For any \( L_1 > 0\) there exists \(L_0 > 0\) such that \(f\) is \((L_0, L_1)\)-smooth.
        \item \label{th:phase-retrieval-smoothness-2} If \(\phi\) satisfies \cref{assumption:phi-star-c2,assumption:isotropic}, then for any \(\bar L > 0\), there exists an \(L > 0\) such that \(f\) satisfies the second-order characterization of \((L, \bar L)\)-anisotropic smoothness.
    \end{theoremenum}
\end{theorem}

\subsubsection{Symmetric matrix factorization}

Consider the symmetric matrix factorization problem with objective and gradient 
\begin{equation} \label{eq:symmetric-factorization-objective-gradient}
    f(U) = \frac{1}{2} \Vert U U^\top - Y \Vert_F^2, \quad \nabla f(U) = (U U^\top - Y) U.
\end{equation}
Here, \(U \in \R^{n \times r}\) is the optimization variable, and \(Y \in \bS^{n \times n}\) is a given symmetric matrix.
When \(r < n\), minimizing \(f\) yields a \emph{low-rank} approximation of \(Y\) with rank at most \(r\).
Such low-rank matrix factorizations are fundamental in a variety of applications, most notably in principal component analysis (PCA) \cite{joliffe_principal_2002}, where one seeks to capture the most significant directions of variation in the data.
More broadly, symmetric matrix factorization plays a central role across various domains: in machine learning, it underlies techniques such as non-negative matrix factorization for parts-based representation learning \cite{lee_learning_1999}; in signal processing, it is employed in matrix completion and compressed sensing to reconstruct structured signals from incomplete or noisy measurements \cite{candes_exact_2009}.

\begin{theorem} \label{th:matrix-factorization-symmetric-smoothness}
    Consider the symmetric matrix factorization problem with objective \eqref{eq:symmetric-factorization-objective-gradient}. Then the following statements hold.
    \begin{theoremenum}
        \item \label{th:matrix-factorization-symmetric-smoothness-1} For any \( L_1 > 0\) there exists an \(L_0 > 0\) such that \(f\) is \((L_0, L_1)\)-smooth.
        \item \label{th:matrix-factorization-symmetric-smoothness-2} If \(\phi\) satisfies \cref{assumption:phi-star-c2,assumption:isotropic}, then for any \(\bar L > 0\), there exists an \(L > 0\) such that \(f\) satisfies the second-order characterization of \((L, \bar L)\)-anisotropic smoothness.
    \end{theoremenum}
\end{theorem}

\subsubsection{Asymmetric matrix factorization}

Consider the regularized asymmetric matrix factorization problem with objective
\begin{equation} \label{eq:asymmetric-factorization-objective}
    f(W, H) = \frac{1}{2} \Vert W H - Y \Vert_F^2 + \frac{\kappa}{4} \Vert W \Vert_F^4 + \frac{\kappa}{4} \Vert H \Vert_F^4,
\end{equation}
where \(W \in \R^{m \times r}\) and \(H \in \R^{r \times n}\) are the optimization variables, \(Y \in \R^{m \times n}\) is a given matrix, and \(\kappa \geq 0\) is a regularization parameter.
When \(\kappa = 0\) and \(r < \min\{m, n\}\), this reduces to the classical low-rank matrix factorization problem.
Additionally, such objectives have been used to model the training of two-layer linear networks, such as in the case of two-layer autoencoders~\cite{fatkhullin_taming_2024}.
We note that the results below also hold for regularization terms of the form \(\kappa \Vert W^\top W - H H^\top \Vert_F^2\) as described in \cite{chi_nonconvex_2019}, and highlight the work of \cite{mukkamala_beyond_2019}, which designed a Bregman proximal-gradient method for similar regularized matrix factorization problems.

\begin{theorem} \label{th:matrix-factorization-asymmetric-smoothness}
    Consider the asymmetric matrix factorization problem with objective \eqref{eq:asymmetric-factorization-objective} and let \(\kappa > 0\). Then the following statements hold.
    \begin{theoremenum}
        \item \label{th:matrix-factorization-asymmetric-smoothness-1} For any \( L_1 > 0\) there exists an \(L_0 > 0\) such that \(f\) is \((L_0, L_1)\)-smooth.
        \item \label{th:matrix-factorization-asymmetric-smoothness-2} If \(\phi\) satisfies \cref{assumption:phi-star-c2,assumption:isotropic}, then for any \(\bar L > 0\), there exists an \(L > 0\) such that \(f\) satisfies the second-order characterization of \((L, \bar L)\)-anisotropic smoothness.
    \end{theoremenum}
\end{theorem}

Note that \cref{th:matrix-factorization-asymmetric-smoothness} requires \(\kappa > 0\). 
To understand why, observe that the gradient of \(f\) is given by
\begin{equation*}
    \nabla_W f(W, H) = (W H - Y) H^\top + \kappa \Vert W \Vert_F^2 W, \quad \text{and} \quad \nabla_H f(W, H) = W^\top (W H - Y) + \kappa \Vert H \Vert_F^2 H.
\end{equation*}
Let \(x\) denote the concatenation of the vectors \(\vecop(W)\) and \(\vecop(H)\), such that \(\Vert x \Vert^2 = \Vert W \Vert_F^2 + \Vert H \Vert_F^2\).
In contrast to symmetric matrix factorization, if \(\kappa = 0\), the gradient norm of \(f\) can approach zero as \( \Vert x \Vert \to \infty \), whereas \cref{assump:novel-condition} requires an asymptotic growth proportional to \(\Vert x \Vert^3\).
To see this, consider \(W^\star \in \R^{m \times r}, H^\star \in \R^{r \times n}\) such that \(W^\star H^\star = Y\). In this case, the gradient norm is zero, and rescaling \(W^\star\) and \(H^\star\) with a nonsingular matrix \(D \in \R^{r \times r}\), i.e., \(\tilde W = W^\star D\) and \(\tilde H = D^{-1} H^\star\), preserves the gradient norm. 
As a result, one can construct counterexamples where the gradient norm remains zero while \(\Vert D \Vert \to \infty\), and consequently \(\Vert x \Vert \to \infty\).

Finally, we remark that the key step in proving \cref{th:matrix-factorization-asymmetric-smoothness} entails lower bounding \(\Vert \nabla f(W, H) \Vert\) in terms of the variable \( V := \max(\Vert W \Vert_F, \Vert H \Vert_F)\), and exploiting that \(\Vert V \Vert \to \infty\) if and only if \(\Vert x \Vert \to \infty\).
It appears that this strategy can be generalized to the factorization of \(Y\) into more than two factors, which is relevant for training deep linear networks.

\subsubsection{Burer-Monteiro factorizations of MaxCut-type semidefinite programs}

Let us consider so-called MaxCut-type semidefinite programs (SDPs)
\begin{equation} \label{eq:maxcut}
    \begin{aligned}
        &\minimize_{X \in \bS^{n \times n}} &&- \langle C, X \rangle\\
        &\stt && X \succeq 0\\
        & &&\diag (X) = 1_n,
    \end{aligned}
\end{equation}
where \(C \in \bS^{n \times n}\) is the \textit{cost matrix}.
The relaxation \eqref{eq:maxcut} provides a precise relaxation to the MaxCut problem, a fundamental combinatorial problem arising in graph optimization \cite{goemans_improved_1995,endor_benign_2025}.
In an effort to exploit the typical low-rank structure of the solution, a \emph{Burer-Monteiro factorization} \cite{burer_nonlinear_2003} decomposes \(X = V V^\top\) for \(V \in \R^{n \times r}\).
This yields
\begin{equation} \label{eq:burer-monteiro}
    \begin{aligned}
        &\minimize_{V \in \bR^{n \times r}} &&- \langle C, V V^\top \rangle\\
        &\stt &&\diag (V V^\top) = 1_n.
    \end{aligned}
\end{equation}
Choosing \(r\) much smaller than \(n\) significantly decreases the number of variables from \(n^2\) to \(nr\).
However, the downside of this approach is that convexity is lost.
Fortunately, under certain conditions every second-order stationary point of this nonconvex problem is a global minimizer \cite{endor_benign_2025}.
Let us denote by \(x_i \in \R^r\) the \(i\)'th row of \(V\), such that \(V^\top = [x_1, x_2, \dots, x_n]\).
We also define the vectorized variable \(x := [x_1^\top, x_2^\top, \dots, x_n^\top]^\top \in \R^{d}\) where \(d = n r\). 
Then we denote by \(f(x)\) the objective of \eqref{eq:burer-monteiro} in terms of \(x\), and likewise by \(A(x) = 0\) the constraint of \eqref{eq:burer-monteiro} in terms of \(x\).

As proposed in the seminal work \cite{burer_nonlinear_2003}, this nonconvex constrained problem can be solved with an augmented Lagrangian method (ALM).
Each iteration consists of minimizing with respect to the primal variable \(x\) the (unconstrained) augmented Lagrangian with penalty parameter \(\beta > 0\), i.e.,
\begin{equation} \label{eq:augmented-lagrangian}
    L_{\beta}(x, y) = f(x) + \langle A(x), y \rangle + \frac{\beta}{2} \Vert A(x) \Vert^2,
\end{equation}
followed by an update of the multipliers \(y \in \R^n\).
A similar strategy was also used in \cite{sahin_inexact_2019} for Burer-Monteiro factorizations of clustering SDPs.
The following theorem establishes generalized smoothness of the augmented Lagrangian with respect to the primal variable.
\begin{theorem}\label{th:burer-monteiro-smoothness}
    Consider the Burer-Monteiro factorization \eqref{eq:burer-monteiro} of the MaxCut-type SDP \eqref{eq:maxcut} and let \(L_{\beta}\) denote the augmented Lagrangian with penalty parameter \(\beta > 0 \) of this factorized problem.
    Then, with respect to the primal variable \(x \in \R^d\) and for some fixed multiplier \(y \in \R^n\) the following statements hold.
    \begin{theoremenum}
        \item \label{th:burer-monteiro-smoothness-1} For any \( L_1 > 0\) there exists \(L_0 > 0\) such that \(L_\beta(\cdot, y)\) is \((L_0, L_1)\)-smooth.
        \item \label{th:burer-monteiro-smoothness-2} If \(\phi\) satisfies \cref{assumption:phi-star-c2,assumption:isotropic}, then for any \(\bar L > 0\), there exists an \(L > 0\) such that \(L_\beta(\cdot, y)\) satisfies the second-order characterization of \((L, \bar L)\)-anisotropic smoothness.
    \end{theoremenum}
\end{theorem}

\section{Saddle point avoidance of the preconditioned gradient method} \label{sec:saddle-point-avoidance}
The remarkable performance of simple gradient descent-like methods for minimizing \emph{nonconvex} functions is often attributed to the fact that they avoid strict saddle points of \textit{Lipschitz smooth} objectives.
This section establishes that nonlinear preconditioning of the gradient preserves this desirable property, and in fact generalizes this result to \emph{anisotropically smooth} functions.

\subsection{Asymptotic results based on the stable-center manifold theorem}

Denote by \( \cX^\star \) the set of strict saddle points of a function \(f \in \cC^2\), i.e., 
\begin{equation*}
    \cX^\star := \left\{ x^\star \mid \nabla f(x^\star) = 0, \quad \lambda_{\min}(\nabla^2 f(x^\star)) < 0 \right\}.
\end{equation*}
Classical results like \cite{lee_gradient_2016,lee_first-order_2019}, which are based on the stable-center manifold theorem \cite{shub_global_1987}, exploit the fact that the eigenvalues of the Hessian \(\nabla^2 f\) are uniformly bounded.
In a similar way, for the preconditioned gradient descent method we require that the second-order characterization of \((L, \bar L)\)-anisotropic smoothness holds.
By exploiting the fact that \( \nabla^2 \phi^*(0) = I \), we then obtain the following theorem, which generalizes \cite[Theorem 4]{lee_gradient_2016}.
\begin{theorem} \label{th:asymptotic-saddle-avoidance-smooth}
    Let \( f \in \cC^2 \) and suppose that \cref{assumption:phi-star-c2} holds.
    Consider the iterates \( \seq{x^k} \) generated by the preconditioned gradient method, i.e., \( x^{k+1} = T_{\gamma, \bar L^{-1}}(x^k) \), where the initial iterate \( x^0 \in \R^n \) is chosen uniformly at random.
    If \(f\) satisfies the second-order characterization of \((L, \bar L)\)-anisotropic smoothness, and if \( \gamma < \frac{1}{L} \), then
    \begin{equation}
        \bP \left( \lim_{k \to \infty} x^k \in \cX^\star \right) = 0.
    \end{equation}
\end{theorem}
\Cref{assumption:phi-star-c2} ensures \(\phi^* \in \cC^2\), which in turn guarantees that \(T_{\gamma, \bar L^{-1}} \in \cC^1\), as needed for the stable-center manifold theorem \cite{shub_global_1987,lee_first-order_2019}.
Unfortunately, this means that the reference function \(\phi(x) = h(\Vert x \Vert)\) with \( h = \frac{1}{2} \Vert \cdot \Vert^2 + \delta_{[-1, 1]} \), which gives rise to a version of the \emph{gradient clipping} method \cite[Example 1.7]{oikonomidis_nonlinearly_2025}, is not covered by \cref{th:asymptotic-saddle-avoidance-smooth}.
Indeed, in this case we have \( {h^*}'(y) = \proj_{[-1, 1]}(y) = \max \left\{ \min \left\{ y, 1 \right\}, -1 \right\} \).
Note however that this projection is a piecewise affine function, and therefore \( {h^*}' \) is continuously differentiable \emph{almost everywhere}, i.e., except at the points \(y = \pm 1\).

Based on a recent variant of the stable-center manifold theorem \cite{cheridito_gradient_2024} we now establish that also the above clipped gradient variant with \(\phi^* \notin \cC^2\) avoids strict saddle points with probability one.
In particular, \cite[Proposition 2.5]{cheridito_gradient_2024} only requires that the iteration map \(T_{\gamma, \lambda}\) is continuously differentiable \emph{on a set of measure one} which contains the set of strict saddle points \(\cX^\star\).
We thus have to show that (i) \(\nabla \phi^*(\bar L^{-1} \nabla f(\cdot)) \) is differentiable almost everywhere; and that (ii) \(\nabla \phi^*(\cdot) \) is differentiable around the point \(\bar L^{-1} \nabla f(x^\star) = 0\), with \(x^\star \in \cX^\star\).
Remark that the former requires an additional assumption for guaranteeing that \(\nabla f\) maps a set of measure one onto a set on which \( \nabla \phi^* \) is differentiable.  

\begin{theorem} \label{th:asymptotic-saddle-avoidance-nonsmooth}
    Let \( f \in \cC^{2+} \) and \(\phi(x) = h(\Vert x \Vert)\) with \( h = \frac{1}{2} \Vert \cdot \Vert^2 + \delta_{[-1, 1]} \).
    Consider the iterates \( \seq{x^k} \) generated by the preconditioned gradient method, i.e., \( 
    x^{k+1} = T_{\gamma, \bar L^{-1}}(x^k) = x^k - \gamma \min (\nicefrac{1}{\Vert \nabla f(x^k) \Vert}, \bar L^{-1}) \nabla f(x^k) 
    \), where the initial iterate \( x^0 \in \R^n \) is chosen uniformly at random.
    Moreover, suppose that the set 
    \begin{equation*}
        U := \left\{ x \in \R^n \mid \Vert \nabla f(x) \Vert \neq \bar L \right\}
    \end{equation*}
    is a set of measure one.
    If \(f\) satisfies the second-order sufficient condition for \((L, \bar L)\)-anisotropic smoothness, and if \( \gamma < \frac{1}{L} \), then
    \begin{equation}
        \bP \left( \lim_{k \to \infty} x^k \in \cX^\star \right) = 0.
    \end{equation}
\end{theorem}

\subsection{Efficiently avoiding strict saddle points through perturbations}

Despite avoiding strict saddle points \emph{asymptotically} for almost any initialization, gradient descent may actually be significantly slowed down around saddle points. 
In fact, gradient descent can take \emph{exponential time} to escape strict saddle points \cite{du_gradient_2017}, in the sense that the number of iterations depends exponentially on the dimension \(n\) of the optimization variable. 
Yet, by adding small perturbations, this issue can be mitigated, and the complexity of obtaining a second-order stationary point then depends only polylogarithmically on the dimension \(n\) \cite{jin_how_2017,jin_nonconvex_2021}.
This section establishes a similar result for a \emph{perturbed preconditioned gradient} method.

Existing works analyzing the \emph{complexity} of gradient descent for converging to a second-order stationary point require not only Lipschitz continuity of the gradients, but also of the Hessian.
This is quite restrictive, since for example any (non-degenerate) polynomial of degree more than \(2\) violates this assumption.
Instead, we require Lipschitz continuity of the mapping
\begin{equation*}
    H_{\lambda}(x) := \lambda^{-1} J [\nabla \phi^*(\lambda \nabla f(x))] = \nabla^2 \phi^*(\lambda \nabla f(x)) \nabla^2 f(x).   
\end{equation*}
To ensure well-definedness of \(H_\lambda\), \cref{assumption:phi-star-c2} is assumed in the remainder of this section.
\begin{assumption} \label{assumption:hessian-lipschitz}
    The mapping \(H_{\lambda}(x) := \nabla^2 \phi^*(\lambda \nabla f(x)) \nabla^2 f(x)\) is \(\rho\)-Lipschitz-continuous, i.e.,
    \begin{equation*}
        \exists \rho >0 : \Vert H_{\lambda}(x) - H_{\lambda}(y) \Vert \leq \rho \Vert x - y \Vert, \quad \forall x, y \in \R^n.
    \end{equation*}
\end{assumption}
This new condition appears significantly less restrictive, as illustrated by the following example.
\begin{example}
    Let \( f(x) = \frac{1}{4} x^4 - \frac{1}{2} x^2 \) and \(\phi(x) = \cosh(\vert x \vert) - 1\).
    Since \(\arcsinh\) is an odd function,
    \begin{equation*}
        \nabla \phi^*(\lambda f'(x)) = \arcsinh(\vert \lambda f'(x) \vert) \barsgn (\lambda f'(x)) = \arcsinh(\lambda f'(x)) = \arcsinh(\lambda (x^3 - x)).
    \end{equation*}
    Therefore, we obtain that \(
        H_{\lambda}(x) = \lambda^{-1} \tfrac{\mathrm d \left( \nabla \phi^*(\lambda f'(x)) \right)}{\mathrm d x} = \tfrac{(3x^2 - 1)}{\sqrt{1+\lambda^2 (x^3-x)^2}}
    \).
    One easily verifies that \(H_{\lambda} \in \cC^1\) with bounded derivative, which implies the required Lipschitz-continuity of \(H_\lambda\).
    In fact, this reasoning generalizes to any univariate polynomial, regardless of its degree.
\end{example}
Under anisotropic smoothness, it is natural to consider \( \lambda^{-1} \phi(\nabla \phi^*(\lambda \nabla f(x))) \) as a first-order stationarity measure, and \(\lambda_{\min}(H_{\lambda}(x))\) as a second-order stationarity measure.
Therefore, we say that a point \(x \in \R^n\) is an \emph{\(\epsilon\)-second-order stationary point} of an \((L, \bar L)\)-anisotropically smooth function \(f\) if
\begin{equation*}
    \lambda^{-1} \phi(\nabla \phi^*(\lambda \nabla f(x))) \leq \epsilon^2, \qquad \text{and} \qquad \lambda_{\min}(\nabla^2 \phi^*(\lambda \nabla f(x)) \nabla^2 f(x)) \geq - \sqrt{\rho \epsilon}.
\end{equation*}
For \(\phi = \frac{1}{2} \Vert \cdot \Vert^2\) we recover the classical notion of \(\epsilon\)-second-order stationarity, with \(\rho\) the constant of Lipschitz continuity of \(\nabla^2 f\).

\Cref{alg:perturbed-GD} describes a perturbed preconditioned gradient method that closely resembles perturbation schemes presented in \cite{jin_how_2017,jin_nonconvex_2021}.
In particular, whenever the first-order stationarity is sufficiently small, then a perturbation is added followed by \(\lceil \cT \rceil > 0\) unperturbed iterations.

\begin{algorithm}
    \caption{Perturbed preconditioned gradient descent}
    \label{alg:perturbed-GD}
    \begin{algorithmic}[1]
        \Require{\(x^0 \in \R^n\), \( \gamma, \lambda > 0\), perturbation radius \(r > 0\), time interval \(\cT >0 \), tolerance \(\cG > 0\)}
        \State \(k_{\text{perturb}} = 0\)
        \For{$k=0, 1, \dots$}
            \If{\(\lambda^{-1} \phi(\nabla \phi^*(\lambda \nabla f(x^k))) \leq \frac{\cG^2}{2}\) and \(k - k_{\text{perturb}} > \cT\)}
                \State \(x^k \leftarrow x^k + \gamma \xi^k\), \quad \(\xi^k \sim \bB_0(r)\) uniformly 
                \State \( k_{\text{perturb}} \leftarrow k\)
            \EndIf
            \State \(x^{k+1} = x^k - \gamma \nabla \phi^*(\lambda \nabla f(x^k))\)
        \EndFor
    \end{algorithmic}
\end{algorithm}\noindent

We analyze the complexity of \cref{alg:perturbed-GD} under the following assumption.
\begin{assumption} \label{assump:perturbed}
    Suppose that \cref{assumption:phi-star-c2} holds, such that \(\phi^* \in \cC^2\), and let \( \phi(x) = h(\Vert x \Vert \)) where in particular \(h \in \cC^2\).
    Moreover, let \(h(x) \geq \nicefrac{x^2}{2}\), and \(h(x) = \nicefrac{x^2}{2} + o(x^2)\) as \( x \to 0\).
\end{assumption}
This assumption holds for kernel functions from \eqref{eq:h-examples}.
Remark that there is no real loss of generality by fixing the scale of \(h\) around \(0\), since a rescaled version of \(h\) can be obtained by modifying \(\bar L\).

In our analysis, we specify the parameters of \cref{alg:perturbed-GD} in terms of \(L, \bar L, \epsilon\) and some \( \chi \geq 1 \),
\begin{equation} \label{eq:perturbed-GD-parameters-1}
    \gamma = \tfrac{1}{L}, \quad \lambda = \tfrac{1}{\bar L}, \quad r = \tfrac{\epsilon}{400 \chi^3}, \quad \cT = \tfrac{L}{\sqrt{\rho \epsilon}} \chi, \quad \cG = \min \left\{ 1, \tfrac{1}{\sqrt{\lambda}} \right\} r,
\end{equation}
and introduce two additional constants that are used only in the analysis, i.e.,
\begin{equation} \label{eq:perturbed-GD-parameters-2}
    \cF = \frac{1}{50\lambda \chi^3} \sqrt{\tfrac{\epsilon^3}{\rho}}, \quad \cZ = \tfrac{1}{4\chi} \sqrt{\tfrac{\epsilon}{\rho}}.
\end{equation}
We obtain the following complexity of \cref{alg:perturbed-GD} for converging to a second-order stationary point.
\begin{theorem}[Iteration complexity] \label{th:saddle-point-complexity}
    Let \(f\) be \((L, \bar L)\)-anisotropically relative to \(\phi\).
    Moreover, suppose that \cref{assump:perturbed,assumption:hessian-lipschitz} hold, and define constants \(\Delta_f \geq f(x^0) - \inf f\) and \( \chi = \log_2 \left(\frac{L^2 \sqrt{n} \Delta_f}{c \sqrt{\rho} \bar L \epsilon^{\nicefrac{5}{2}} \delta} \right)\) for some \( c > 0\). 
    There exists a constant \(c_{\max} > 0\) such that if \( c \leq c_{\max}\), then for any \( \epsilon > 0\) sufficiently small, and for any \( \delta \in (0, 1)\), \Cref{alg:perturbed-GD} with parameters as in \eqref{eq:perturbed-GD-parameters-1} and \eqref{eq:perturbed-GD-parameters-2}, visits an \(\epsilon\)-second-order stationary point in at least \( \nicefrac{T}{2}\) iterations with probability at least \(1 - \delta\), where
    \begin{equation*}
        T = 8 \max \left\{
            \frac{(f(x^0) - \inf f) \cT}{\cF}, \lambda \frac{(f(x^0) - \inf f)}{2 \gamma \cG^2}
        \right\} = \tilde O \left( \frac{L (f(x^0) - \inf f)}{\bar L \epsilon^2} \right).
    \end{equation*}
\end{theorem}
The \(\tilde \cO\) notation hides a factor \(\chi^4\) which is polylogarithmic in the dimension \(n\) and in the tolerance \(\epsilon\).

\Cref{th:saddle-point-complexity} generalizes \cite[Theorem 18]{jin_nonconvex_2021}, and relies on a similar high-level proof strategy, which goes as follows.
If the current iterate \( x \) is not an \(\epsilon\)-second-order stationary point, then either \( \lambda^{-1} \phi(\nabla \phi^*(\lambda \nabla f(x))) \) is large, or \(\lambda_{\min}(H_{\lambda})\) is sufficiently negative.
In either case, we establish a significant decrease in function value after at most \(\lceil \cT \rceil\) iterations of \cref{alg:perturbed-GD}.
Since \( f(x^0) - \inf f \) is bounded, the number of iterates which are not \(\epsilon\)-second-order stationary can be bounded.

Nevertheless, the generalization of \cite[Theorem 18]{jin_nonconvex_2021} to the setting of \Cref{alg:perturbed-GD} is by no means straightforward.
The original proofs rely heavily on Lipschitz smoothness, in a way that often does not generalize directly to the anisotropically smooth setting.
Here, we highlight two such difficulties.
First, consider a point \(x \in \R^n\) and the perturbed point \(\bar x := x + \gamma \xi\) for some perturbation \(\xi \in \bB_0(r)\).
Then, by anisotropic smoothness we can upper bound 
\begin{equation*}
    f(\bar x) - f(x) \leq \frac{\gamma}{\lambda} \phi(\xi + \nabla \phi^*(\lambda \nabla f(x))).
\end{equation*}
While Lipschitz-smoothness with \(\phi = \frac{1}{2} \Vert \cdot \Vert^2\) readily provides an upper bound in terms of \( \Vert \xi \Vert^2 \leq r^2 \) and \(\Vert \nabla f(x) \Vert^2\), the reference functions are not typically such that an upper bound in terms of
\begin{equation*}
    \phi(\xi) + \phi (\nabla \phi^*(\lambda \nabla f(x))) 
\end{equation*}
can be obtained.
Second, unlike in the \(L\)-Lipschitz-smooth case where \( \Vert \nabla^2 f \Vert \leq L\), the norm \(\Vert H_{\lambda} \Vert\) cannot be upper bounded uniformly, even under the second-order characterization of \((L, \bar L)\)-anisotropic smoothness. 
The latter only guarantees that \(\lambda_{\max}(H_\lambda) \leq L \bar L \), but it does not lower bound \(\lambda_{\min}(H_\lambda)\).
And even if the eigenvalues of \(H_\lambda\) were bounded in absolute value, this still would not guarantee boundedness of \(\Vert H_\lambda \Vert\), since \(H_\lambda\) is not a normal matrix in general.

\section{Numerical validation}
\begin{figure}
        \centering
        \includegraphics[width=0.96\textwidth]{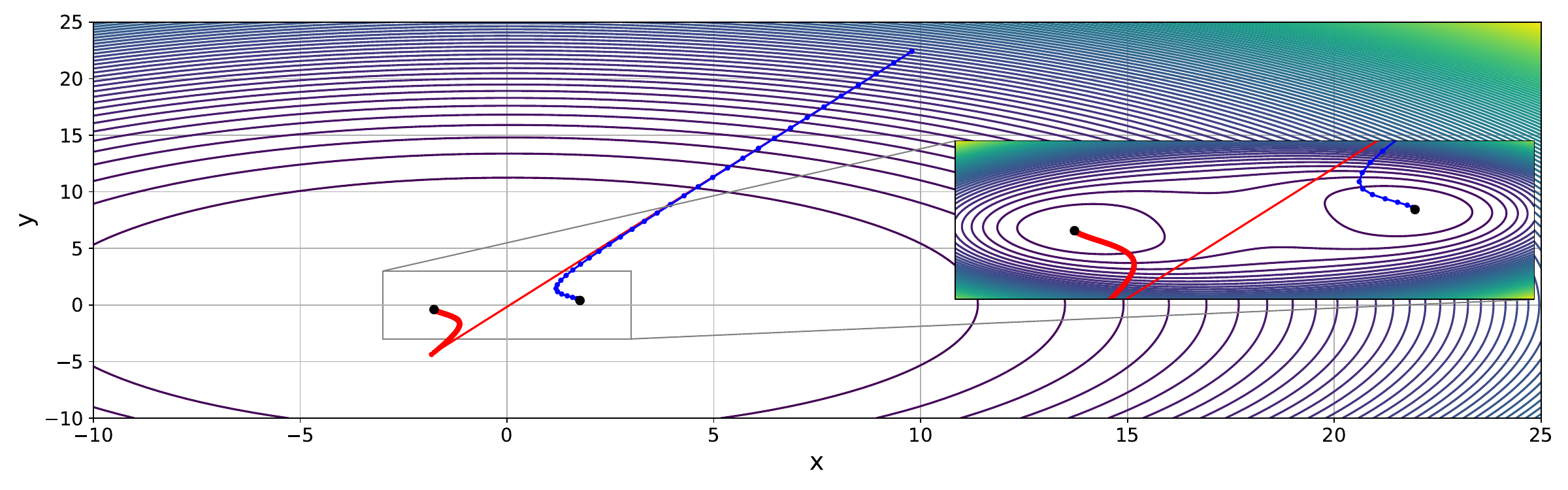}
        \caption{
        Iterates of GD (red) and \eqref{eq:pgd} (blue) on a symmetric matrix factorization problem.
        }
        \label{fig:matrix-factorization}
    \end{figure}

\paragraph{Nonlinear preconditioning for symmetric matrix factorization}
For the symmetric matrix factorization problem \eqref{eq:symmetric-factorization-objective-gradient} with $n = 2$ and $r = 1$, Fig.\;\ref{fig:matrix-factorization} presents a 2D visualization of the level curves of the objective, along with the iterates of both vanilla gradient descent (GD) and the preconditioned variant \eqref{eq:pgd} with $\phi(x) = \cosh(\Vert x \Vert) - 1$. Unless GD is initialized close to a stationary point, the stepsize must be chosen very small to prevent the iterates from diverging -- as expected, because the quartic objective is not Lipschitz smooth. In contrast, the \eqref{eq:pgd} iterations take the form (for $\bar L = 1$)
$$
x^+ = x - \gamma \frac{\sinh^{-1}(\Vert \nabla f(x) \Vert)}{\Vert \nabla f(x) \Vert} \nabla f(x).
$$
In this case, \emph{large gradients are damped} -- recall the close resemblance to clipping methods, cf. Fig.\;\ref{fig:kernels} -- resulting in the convergence of \eqref{eq:pgd} for stepsizes $\gamma$ that are often \emph{orders of magnitudes} larger than the maximum stepsize of GD. In turn, this causes \eqref{eq:pgd} to often require significantly fewer iterations, and overall outperform GD for fixed stepsize.

\begin{figure}[t]
    \centering
        \includegraphics[width=0.24\linewidth]{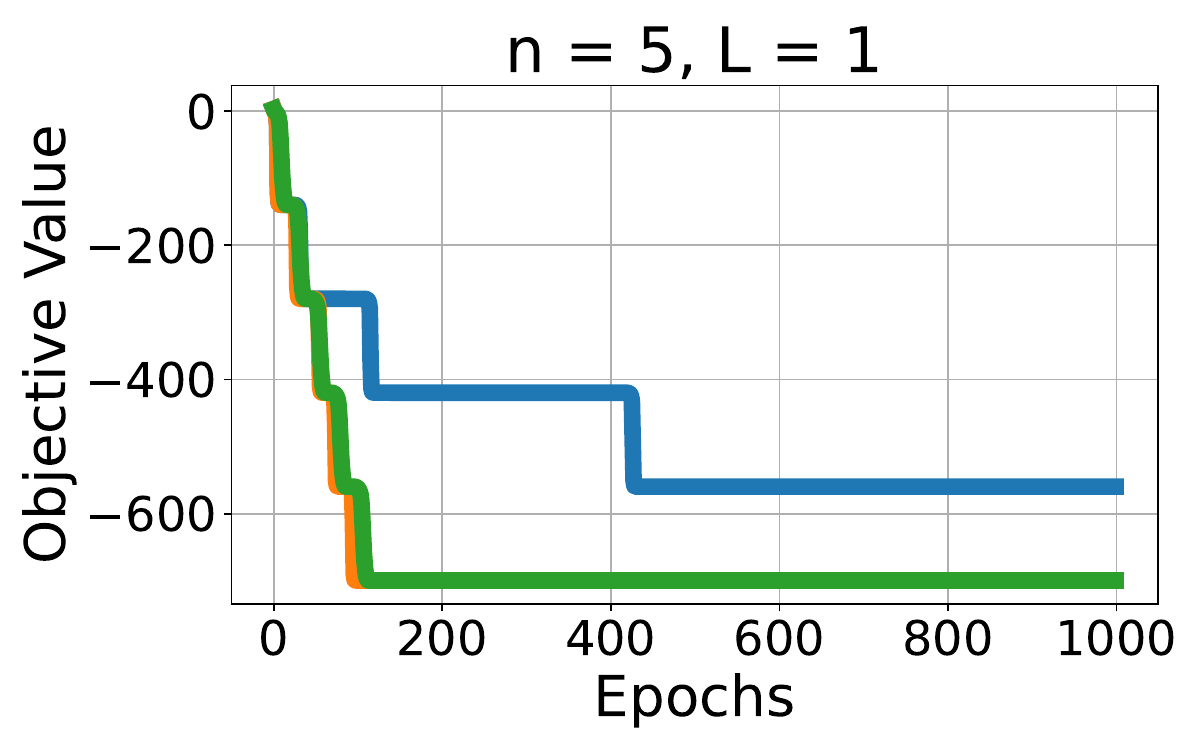}
    \hfill
        \includegraphics[width=0.24\linewidth]{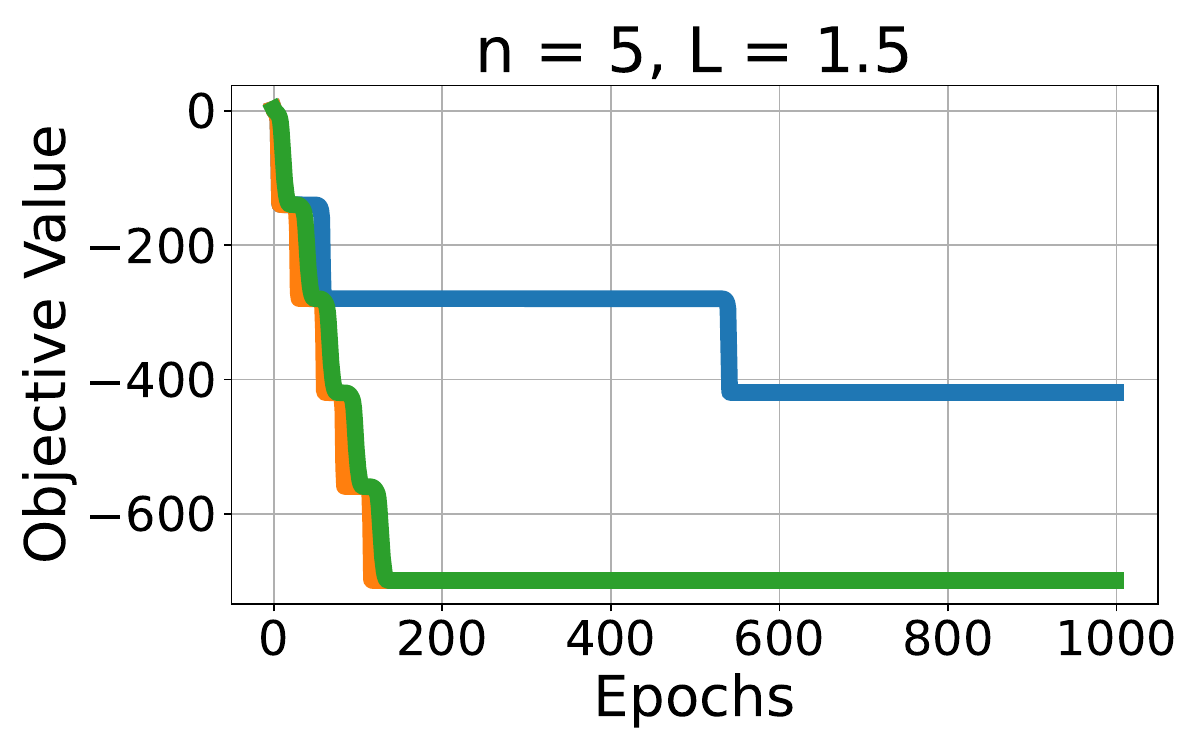}
    \hfill
        \includegraphics[width=0.24\linewidth]{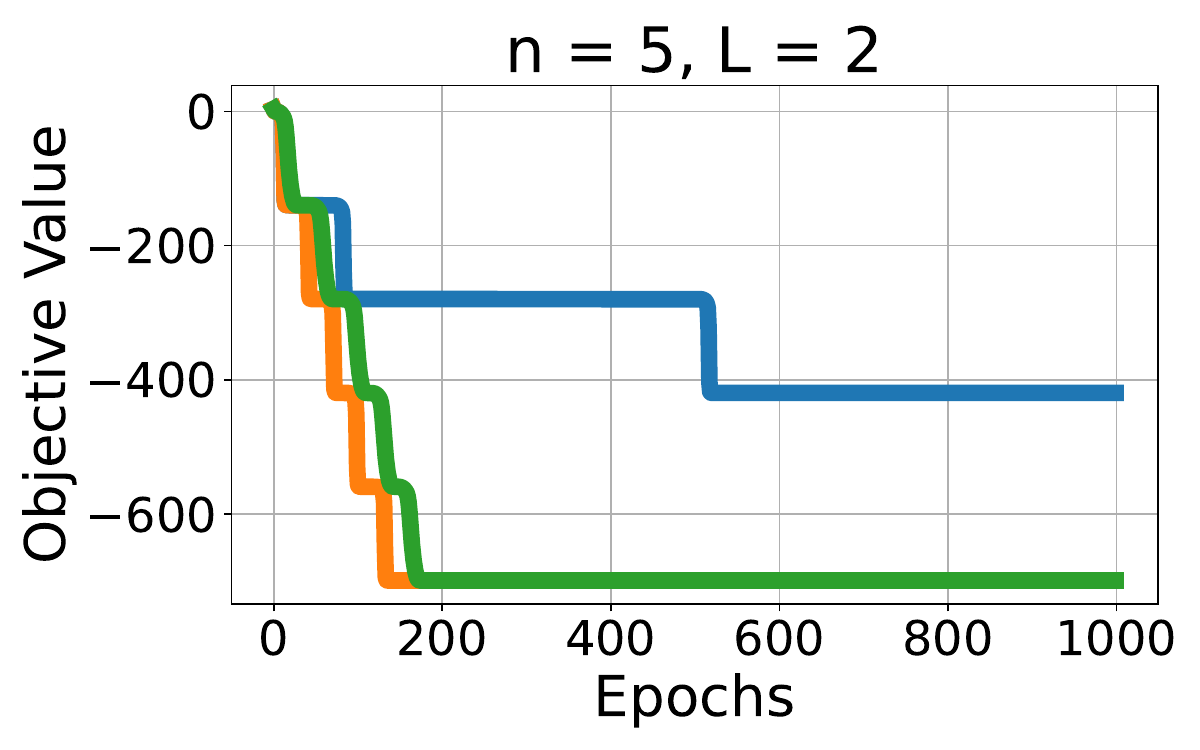}
    \hfill
        \includegraphics[width=0.24\linewidth]{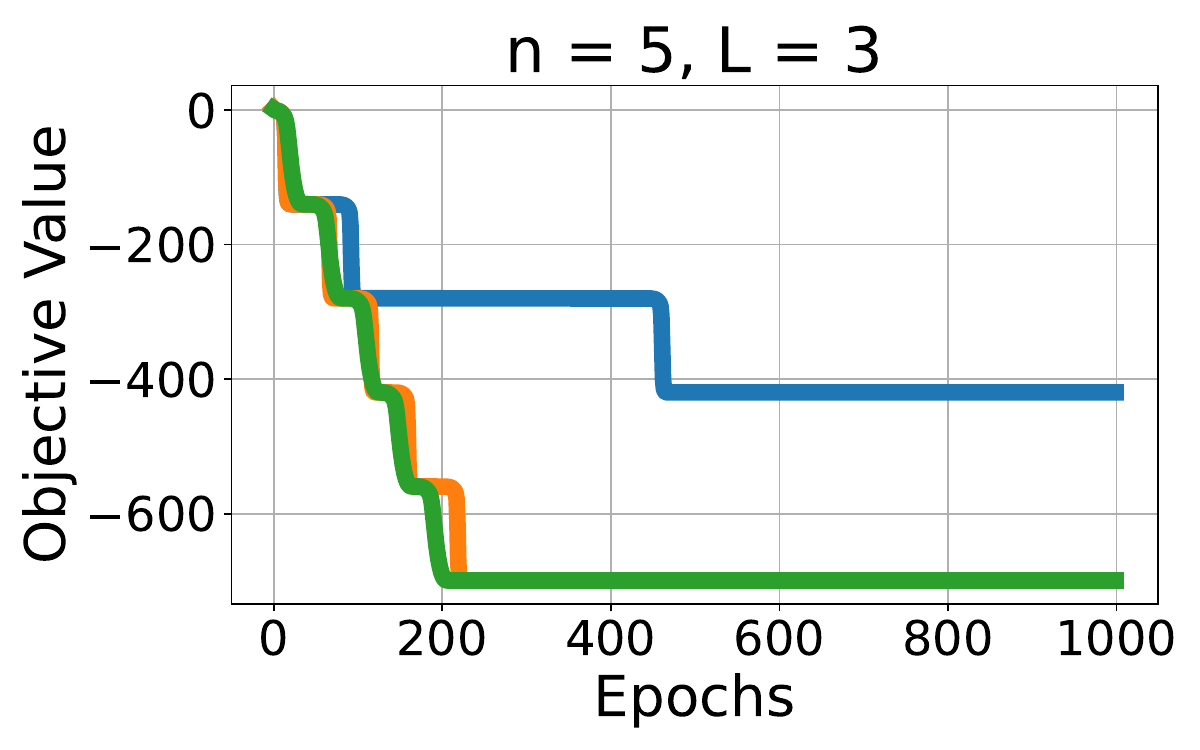}
    \hfill
        \includegraphics[width=0.24\linewidth]{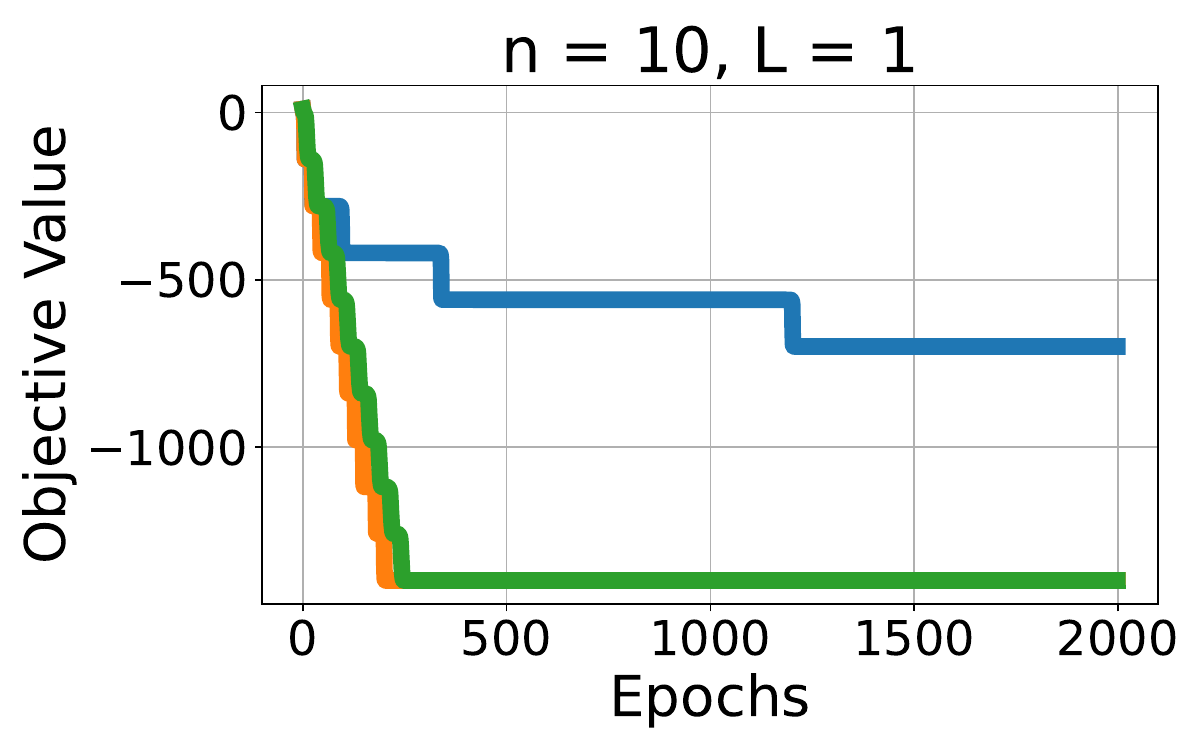}
    \hfill
        \includegraphics[width=0.24\linewidth]{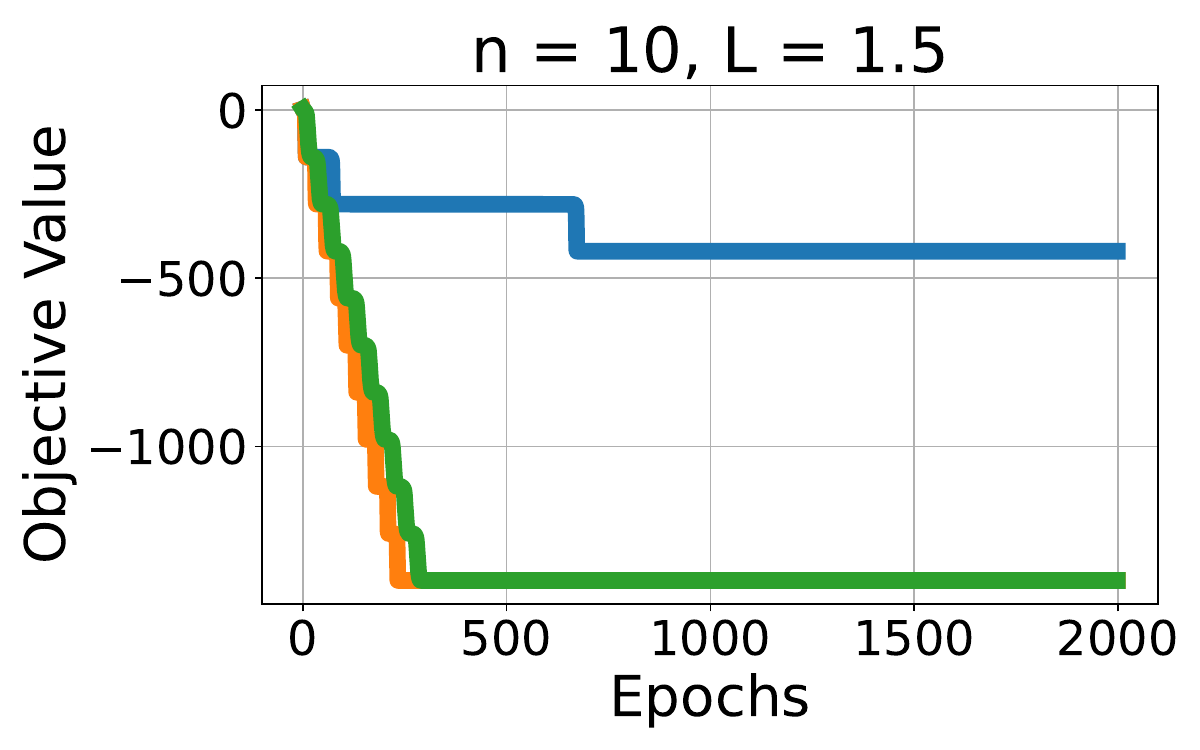}
    \hfill
        \includegraphics[width=0.24\linewidth]{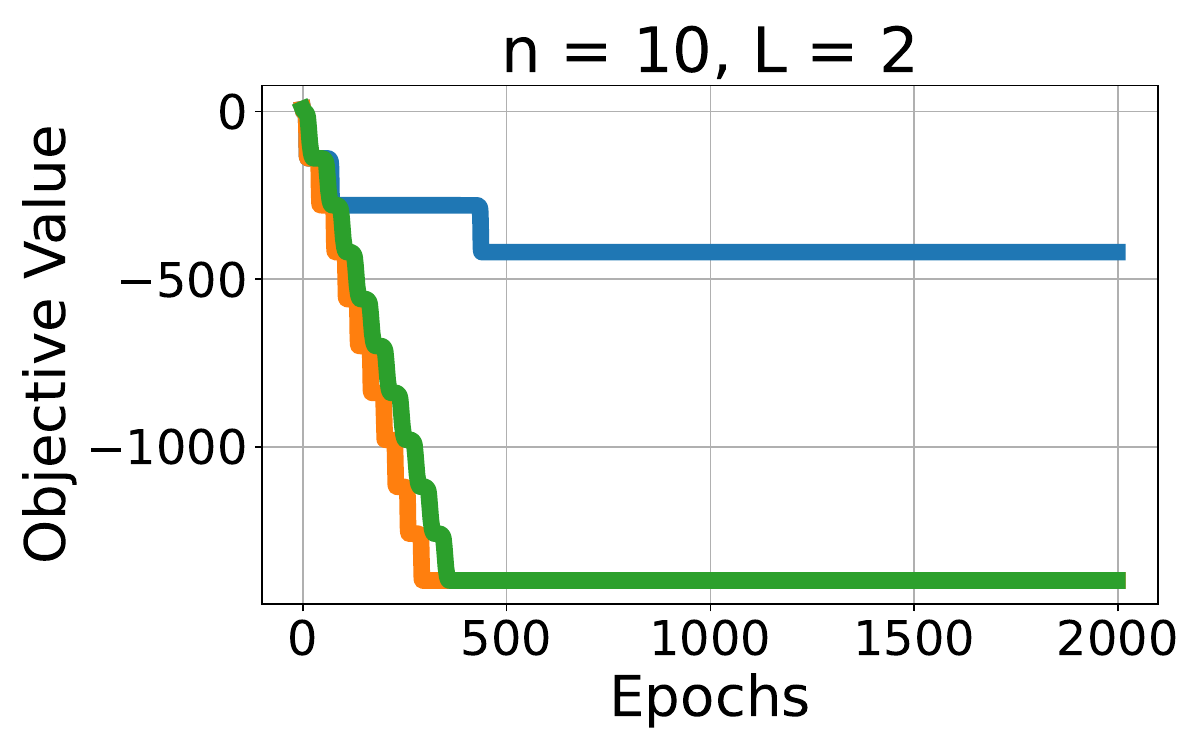}
    \hfill
        \includegraphics[width=0.24\linewidth]{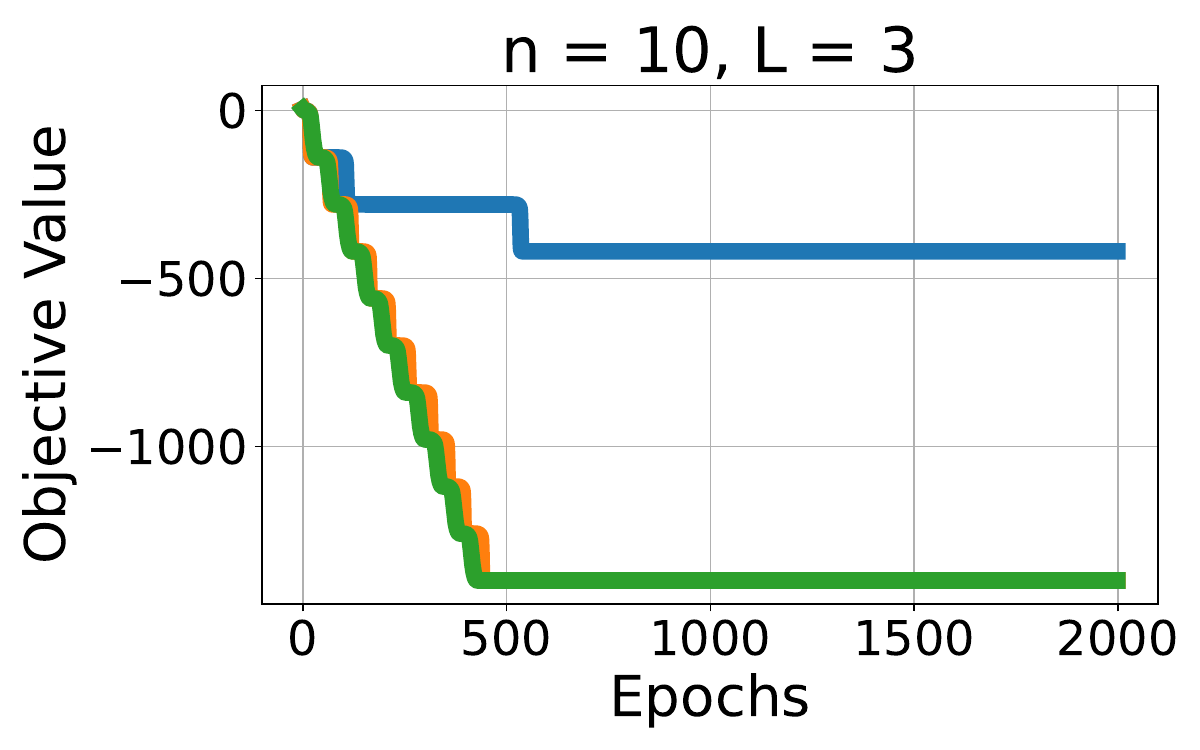}
    \caption{Performance of vanilla GD (blue), perturbed vanilla GD \cite[Alg\;1]{du_gradient_2017} (orange), and \Cref{alg:perturbed-GD} (green) on the `octopus' function \cite{du_gradient_2017}.}
    \label{fig:octopus}
\end{figure}

\paragraph{Fast avoidance of saddle points} 
Fig\;\ref{fig:octopus} validates the fast escape of saddle points by \Cref{alg:perturbed-GD}. We consider the `octopus' objective \cite{du_gradient_2017} which was constructed such that GD takes exponential time to escape saddle points. We select all hyperparameters as in \cite[\S 5]{du_gradient_2017}, and set the only additional hyperparameter $\bar L = 1$. We compare against vanilla GD and perturbed vanilla GD \cite[Alg\;1]{du_gradient_2017}, and vary the constant $L \in \{1, 1.5, 2, 3\}$ and dimension $n \in \{5, 10\}$, thus creating counterparts to [11, Figs 3 and 4].
We observe that \cref{alg:perturbed-GD} performs similar to perturbed vanilla GD, and also scales in a similar way with respect to \(n\) and \(L\). This validates the complexity result from \cref{th:saddle-point-complexity}.

\section{Conclusion}

This work introduced a novel sufficient condition unifying \((L_0, L_1)\)-smoothness and anisotropic smoothness. We showed that this condition holds in key applications such as phase retrieval, matrix factorization, and Burer-Monteiro factorizations of MaxCut.

We further analyzed the nonlinearly preconditioned gradient method, which naturally aligns with anisotropic smoothness. Notably, we proved that it preserves the saddle point avoidance properties of gradient descent and extends them to anisotropically smooth settings. This contrasts with prior analyses requiring either global Lipschitz smoothness, or local smoothness combined with compactness, both of which are often unmet in practice.

To our knowledge, this is the first work to rigorously establish saddle point avoidance for problems like phase retrieval and matrix factorization under a smoothness condition that is both practical and verifiable. These results strengthen the theoretical foundations of first-order methods for nonconvex optimization and in particular encourage further study of nonlinear gradient preconditioning.

\begin{ack}
    The authors are supported by the Research Foundation Flanders (FWO) research projects G081222N, G033822N, G0A0920N and the Research Council KU Leuven C1 project with ID C14/24/103.
\end{ack}

\printbibliography

\newpage
\appendix
\section{Additional results} 
\subsection{Univariate polynomials satisfy Assumption \ref{assump:novel-condition}} \label{sec:univariate-polynomials-verification}

\begin{theorem}
    Let \(f(x) = \sum_{i = 0}^d a_i x^i\) be a univariate polynomial of degree \(d\) in \(x\) with coefficients \(a_i \in \R\).
    Then, \(f\) satisfies \cref{assump:novel-condition}.
\end{theorem}
\begin{proof}
    Without loss of generality, assume that \(a_d \neq 0\), since otherwise \(f\) would be a polynomial of lower degree. 
    Then, by the triangle inequality we have
    \begin{equation*}
        \vert f''(x) \vert \leq \sum_{i = 0}^{d-2} \vert (i+1) (i+2) a_{i+2} x^i \vert \leq \sum_{i = 0}^{d-2} (i+1) (i+2) \vert a_{i+2} \vert \vert x \vert^i
    \end{equation*}
    which is a polynomial of degree \(d-2\) in \(\vert x \vert\). In a similar way, we obtain from the triangle inequality 
    \begin{equation*}
        \begin{aligned}
            \vert f'(x) \vert \geq \vert d a_{d} x^{d-1} \vert - \vert \sum_{i = 0}^{d-2} (i+1) a_{i+1} x^i \vert &\geq \vert d a_{d} x^{d-1} \vert - (d-1)\sum_{i = 0}^{d-2} \vert a_{i+1} \vert \vert x^i \vert\\
            &= d \vert a_d \vert \vert x^{d-1} \vert - (d-1) \sum_{i = 0}^{d-2} \vert a_{i+1} \vert \vert x^i \vert
        \end{aligned}
    \end{equation*}
    which is a polynomial of degree \( d - 1\) in \(\vert x \vert\) where the leading coefficient \(d \vert a_d \vert\) is nonzero.
\end{proof}

\subsection{Multivariate polynomials for which \texorpdfstring{\((L_0, L_1)\)-smoothness}{(L0, L1)-smoothness} fails} \label{sec:multivariate-polynomials-fail-extended}

This section extends \cref{ex:multivariate} and provides some simple multivariate polynomials which are not \((L_0, L_1)\)-smooth.
In particular, we illustrate that this may still happen if the gradient norm grows unbounded.

\begin{figure}[ht]
    \centering
    \includegraphics[width=\linewidth]{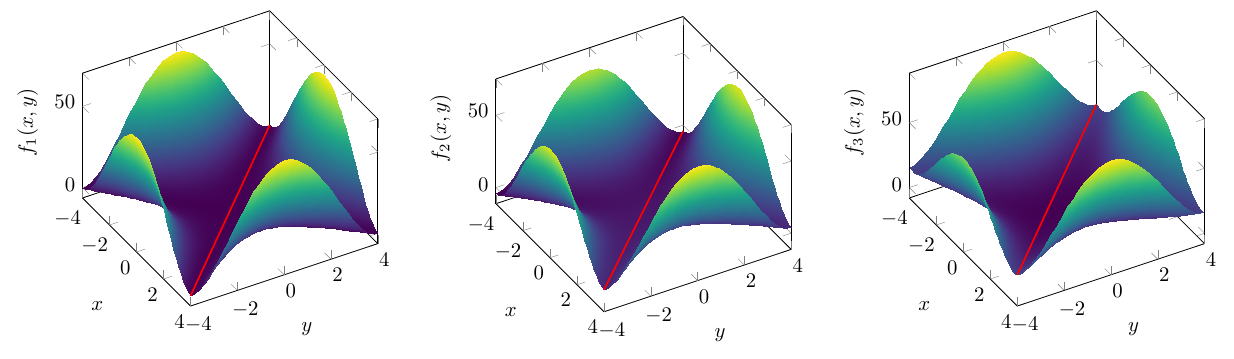}
    \caption{Surface plot of some multivariate polynomials which are not \((L_0, L_1)\)-smooth. The gradient norm is zero (left), constant (middle), or scales proportional to \(\vert x \vert\) (right) along the path \(y = -x\) (red), whereas the Hessian norm scales with \(\vert x \vert^2\).}
    \label{fig:surface_plot}
\end{figure}

Consider the following functions
\begin{equation*}
    f_1(x, y) = \frac{1}{4}(x^4+y^4) - \frac{1}{2} x^2 y^2, \qquad f_2(x, y) = f_1(x, y) + x, \qquad f_3(x, y) = f_1(x, y) + x^2.
\end{equation*}
By a similar reasoning as in \cref{ex:multivariate}, we remark that along a path \(y = -x\) these functions have gradients
\begin{equation*}
    \nabla f_1(x, -x) = \begin{pmatrix}
        0\\
        0
    \end{pmatrix}, \qquad \nabla f_2(x, -x) = \begin{pmatrix}
        1\\
        0
    \end{pmatrix}, \qquad \nabla f_3(x, -x) = \begin{pmatrix}
        2x\\
        0
    \end{pmatrix}
\end{equation*}
and Hessians
\begin{equation*}
    \nabla^2 f_1(x, -x) = \nabla^2 f_2(x, -x) = x^2 \begin{pmatrix}
        2 & -2\\
        -2 & 2
    \end{pmatrix}, \quad \nabla^2 f_3(x, -x) = x^2 \begin{pmatrix}
        2 & -2\\
        -2 & 2
    \end{pmatrix} + 2 \begin{pmatrix}
        1 & 0\\
        0 & 0
    \end{pmatrix}.
\end{equation*}
Clearly, these functions cannot be \((L_0, L_1)\)-smooth, because for \(y = -x\) the Hessian norms grow proportionally to \(\vert x \vert^2\), whereas the gradient norms are zero (for \(f_1\)), constant (for \(f_2\)), or grow proportionally to \(\vert x \vert\) (for \(f_3\)).
This is visualized in \cref{fig:surface_plot}.

Remark that \(f_3\) illustrates that unboundedness of the gradient norm is not sufficient for \((L_0, L_1)\)-smoothness.
Instead, the gradient norm needs to grow `sufficiently fast'; a sufficient condition is given by \cref{assump:novel-condition}.

\subsection{Anisotropic smoothness is more general than \texorpdfstring{\((L_0, L_1)\)-smoothness}{(L0, L1)-smoothness}}

\Cref{ex:l0-l1-smoothness} already established that if a function \(f\) is \((L_0, L_1)\)-smooth, then it also satisfies the second-order characterization of \((L_1, \nicefrac{L_0}{L_1})\)-anisotropic smoothness (\cref{assump:novel-condition}) relative to the reference function \(\phi(x) = - \Vert x \Vert - \ln (1 - \Vert x \Vert)\).
We now show that the function \(f(x) = \exp(\Vert x \Vert^2)\) is \emph{not \((L_0, L_1)\)-smooth}, but \emph{satisfies the second-order characterization of \((L, \bar L)\)-anisotropic smoothness} with \(L = 2, \bar L = 1\), thus confirming that anisotropic smoothness generalizes \((L_0, L_1)\)-smoothness.

Note that the function \(f(x) = \exp(\Vert x \Vert^2)\) has gradient and Hessian
\begin{equation*}
    \nabla f(x) = 2 \exp(r^2) x, \qquad \nabla^2 f(x) = 2 \exp(r^2) \left( I + 2 x x^\top \right)
\end{equation*}
where we defined \(r := \Vert x \Vert\) for ease of notation.
Remark also that 
\begin{equation*}
    \Vert \nabla f(x) \Vert = 2 \exp(r^2) \Vert x \Vert = 2 r \exp(r^2), \qquad \Vert \nabla^2 f(x) \Vert = 2 \exp(r^2)(1+2r^2),
\end{equation*}
where the norm of the Hessian follows from the observation that \(x x^\top\) has eigenvalues \(0\) and \(r^2\).

The following theorem establishes that \(f(x) = \exp(\Vert x \Vert^2)\) is not \((L_0, L_1)\)-smooth.
\begin{theorem}
    There do not exist constants \(L_0, L_1 \geq 0\) for which the function \(f(x) = \exp(\Vert x \Vert^2)\) is \((L_0, L_1)\)-smooth.
\end{theorem}
\begin{proof}
    Assume, by contradiction that there exist constants \(L_0, L_1 \geq 0\) for which
    \begin{equation*}
        \Vert \nabla^2 f(x) \Vert \leq L_0 + L_1 \Vert \nabla f(x) \Vert, \qquad \text{for all } x \in \R^n.
    \end{equation*}
    This means that 
    \begin{equation*}
        2 \exp(r^2)(1+2r^2) \leq L_0 + 2 r L_1 \exp(r^2), \qquad \text{for all } x \in \R^n,
    \end{equation*}
    or equivalently,
    \begin{equation*}
        1 + 2 r^2 \leq \frac{L_0}{2 \exp(r^2)} +L_1 r, \qquad \text{for all } x \in \R^n.
    \end{equation*}
    Clearly, this cannot hold, since the left hand side grows faster than the right hand side as \(r \to \infty\).
\end{proof}

\begin{theorem}
    The function \(f(x) = \exp(\Vert x \Vert^2)\) satisfies the second-order characterization of \((L, \bar L)\)-anisotropic smoothness relative to the reference function \( \phi(x) = -\Vert x \Vert - \ln(1 - \Vert x \Vert) \) for \(\bar L = 1\) and \(L \geq 2\).
\end{theorem}
\begin{proof}
    By \cite[Lemma 2.5]{oikonomidis_nonlinearly_2025}, the second-order characterization of anisotropic smoothness for \(\bar L = 1\) is equivalent to
    \begin{equation*}
        L^{-1} \nabla^2 f(x) \preceq \left[ \nabla^2 \phi^*(\nabla f(x)) \right]^{-1}.
    \end{equation*}
    For our particular reference function, this condition becomes
    \begin{equation*}
        L^{-1} \nabla^2 f(x) \preceq (1+\Vert \nabla f(x) \Vert) (I + \tfrac{1}{\Vert \nabla f(x) \Vert}\nabla f(x) \nabla f(x)^\top). 
    \end{equation*}
    Let us define \(\alpha := 1+ \Vert \nabla f(x) \Vert = 1+ 2r \exp(r^2)\) and recall that \(\nabla f(x) = 2 \exp(r^2) x\).
    Then it remains to show that the matrix 
    \begin{equation*}
        M_L(x) = \alpha \left( I + \Vert \nabla f(x) \Vert \frac{x x^\top}{ \Vert x \Vert^2} \right) - \frac{2 \exp(r^2)}{L} \left(I + 2 x x^\top \right)
    \end{equation*}
    is positive semidefinite for \(L \geq 2\), uniformly in \(x \in \R^n\).
    Note that we can rewrite
    \begin{equation*}
        M_L(x) = \underbrace{\left( \alpha - \frac{2 \exp(r^2)}{L} \right)}_{:= A} I + \underbrace{\left( \alpha \Vert \nabla f(x) \Vert - \frac{4 \exp(r^2)}{L} \Vert x \Vert^2 \right)}_{:= B} \frac{x x^\top}{ \Vert x \Vert^2}
    \end{equation*}
    as the weighted sum of two symmetric positive semidefinite matrices.
    We now show that both weights \(A\) and \(B\) are nonnegative for \(L \geq 2\), for any \(r \geq 0\), from which the claim follows.

    The weight \(A\) equals
    \begin{equation*}
        A(r) = \alpha - \frac{2 \exp(r^2)}{L} = 1+ \exp(r^2) \left( 2r - \frac{2}{L} \right)
    \end{equation*}
    and has derivative
    \begin{equation*}
        A'(r) = \exp(r^2) \left( 2 - \frac{4r}{L} + 4 r^2 \right) = \exp(r^2) \left( 2 - \frac{1}{L^2} + \left( \frac{1}{L} -  2 r \right)^2 \right).
    \end{equation*}
    For \(L \geq 2\) this yields
    \begin{equation*}
        A'(r) \geq \exp(r^2) \left( 1 + \left( \frac{1}{L} -  2 r \right)^2 \right) > 0,
    \end{equation*}
    meaning that \(A(r)\) is strictly increasing. 
    Since \(A(0) = 1 - \frac{2}{L}\), we conclude that \(A(r)\) is nonnegative for \(r \geq 0\).

    As for the weight \(B\), remark that \[
        \alpha \Vert \nabla f(x) \Vert = (1+2r \exp(r^2)) 2r \exp(r^2) = 2r \exp(r^2) + 4 r^2 \exp(2r^2),
    \]
    and \[
        \frac{4 \exp(r^2)}{L} \Vert x \Vert^2 = \frac{4 r^2 \exp(r^2)}{L}.
    \]
    Therefore, the second coefficient becomes
    \begin{equation*}
        B(r) = 4 r^2 \exp(r^2) \left( \exp(r^2) - \frac{1}{L} \right) + 2r \exp(r^2)
    \end{equation*}
    which is nonnegative for \(L \geq 2\) since \(\exp(r^2) \geq 1 > \frac{1}{L}\).
\end{proof}

As a second example, we consider the function \(f(x) = \exp(\Vert x \Vert^2) - 2 \Vert x \Vert^2\), which has gradient and Hessian
\begin{equation*}
    \nabla f(x) = 2 \left( \exp(r^2) - 2 \right) x, \qquad \nabla^2 f(x) = 2 \exp(r^2) \left( I + 2 x x^\top \right) - 4 I
\end{equation*}
where we defined \(r := \Vert x \Vert\) for ease of notation.
Remark also that 
\begin{equation*}
    \Vert \nabla f(x) \Vert = 2 \left\vert \exp(r^2) - 2 \right\vert \Vert x \Vert = 2 r \left\vert \exp(r^2) - 2 \right\vert, \quad \Vert \nabla^2 f(x) \Vert = 2 \exp(r^2)(1+2r^2) - 4,
\end{equation*}
where the norm of the Hessian follows from the observation that \(x x^\top\) has eigenvalues \(0\) and \(r^2\).

The following theorem establishes that \(f\) is not \((L_0, L_1)\)-smooth.
\begin{theorem}
    There do not exist constants \(L_0, L_1 \geq 0\) for which the function \(f(x) = \exp(\Vert x \Vert^2) - 2 \Vert x \Vert^2\) is \((L_0, L_1)\)-smooth.
\end{theorem}
\begin{proof}
    Assume, by contradiction that there exist constants \(L_0, L_1 \geq 0\) for which
    \begin{equation*}
        \Vert \nabla^2 f(x) \Vert \leq L_0 + L_1 \Vert \nabla f(x) \Vert, \qquad \text{for all } x \in \R^n.
    \end{equation*}
    This means that 
    \begin{equation*}
        2 \exp(r^2)(1+2r^2) - 4 \leq L_0 + 2 r L_1 \left\vert \exp(r^2) - 2 \right\vert, \qquad \text{for all } x \in \R^n,
    \end{equation*}
    or equivalently,
    \begin{equation*}
        1 + 2 r^2 \leq \frac{L_0 + 4}{2 \exp(r^2)} +L_1 r \frac{\left\vert \exp(r^2) - 2 \right\vert}{\exp(r^2)}, \qquad \text{for all } x \in \R^n.
    \end{equation*}
    Clearly, this cannot hold, since the left hand side grows faster than the right hand side as \(r \to \infty\).
\end{proof}

\begin{theorem}
    The function \(f(x) = \exp(\Vert x \Vert^2) - 2 \Vert x \Vert^2\) satisfies the second-order characterization of \((L, \bar L)\)-anisotropic smoothness relative to the reference function \( \phi(x) = -\Vert x \Vert - \ln(1 - \Vert x \Vert) \) for \(\bar L = 1\) and \(L = 10\).
\end{theorem}
\begin{proof}
    By \cite[Lemma 2.5]{oikonomidis_nonlinearly_2025}, the second-order characterization of anisotropic smoothness for \(\bar L = 1\) is equivalent to
    \begin{equation*}
        L^{-1} \nabla^2 f(x) \preceq \left[ \nabla^2 \phi^*(\nabla f(x)) \right]^{-1}.
    \end{equation*}
    For our particular reference function, this condition becomes
    \begin{equation*}
        L^{-1} \nabla^2 f(x) \preceq (1+\Vert \nabla f(x) \Vert) (I + \tfrac{1}{\Vert \nabla f(x) \Vert} \nabla f(x) \nabla f(x)^\top). 
    \end{equation*}
    Let us define \(\alpha := 1+ \Vert \nabla f(x) \Vert = 1+ 2 r \left\vert \exp(r^2) - 2 \right\vert\) and recall that \(\nabla f(x) = 2 \left( \exp(r^2) - 2 \right) x\).
    Then it remains to show that the matrix 
    \begin{equation*}
        M_L(x) = \alpha \left( I + \Vert \nabla f(x) \Vert \frac{x x^\top}{ \Vert x \Vert^2} \right) - \frac{2 \exp(r^2)}{L} \left(I + 2 x x^\top \right) + \frac{4}{L} I
    \end{equation*}
    is positive semidefinite for \(L = 10\), uniformly in \(x \in \R^n\).
    Note that we can rewrite
    \begin{equation*}
        M_L(x) = \underbrace{\left( \alpha - \frac{2 \exp(r^2)}{L} + \frac{4}{L} \right)}_{:= A} I + \underbrace{\left( \alpha \Vert \nabla f(x) \Vert - \frac{4 \exp(r^2)}{L} \Vert x \Vert^2 \right)}_{:= B} \frac{x x^\top}{ \Vert x \Vert^2}.
    \end{equation*}
    Our proof strategy goes as follows.
    First, we show that \(A(r) \geq 0.8\) for all \(r \geq 0\).
    Since the eigenvalues of \(\frac{v v^\top}{ \Vert v \Vert^2}\) are \(0\) and \(1\), and since adding a multiple of the identity matrix shifts the eigenvalues by that multiple, it then suffices to show that \(B(r) \geq -0.8\) for \(r \geq 0\).
    This then implies positive semidefiniteness of \(M_L\) and proves the claim.

    We start by lower bounding
    \begin{equation*}
        A(r) = \alpha - \frac{2 \exp(r^2)}{L} + \frac{4}{L} = 1 + 2 r \left\vert \exp(r^2) - 2 \right\vert - \frac{2 \exp(r^2)}{L} + \frac{4}{L}.
    \end{equation*}
    We distinguish three cases.
    If \(\exp(r^2) < 2\), then \(\left\vert \exp(r^2) - 2 \right\vert = 2 - \exp(r^2)\), and hence
    \begin{align*}
        A(r) &= 1 + 2 r \left( 2 - \exp(r^2) \right) - \tfrac{2 \exp(r^2)}{L} + \frac{4}{L} = 1 + 4r - 2 \exp(r^2) \left(r + \frac{1}{L}\right) + \frac{4}{L}\\
        &> 1 + 4 r - 4 \left(r + \frac{1}{L} \right) + \frac{4}{L} = 1 - \frac{4}{L} + \frac{4}{L} = 1.
    \end{align*}
    On the other hand, if \(\exp(r^2) \geq 2\), then \(\left\vert \exp(r^2) - 2 \right\vert = \exp(r^2) - 2\), and hence
    \begin{align*}
        A(r) &= 1 + 2 r \left( \exp(r^2) - 2 \right) - \tfrac{2 \exp(r^2)}{L} + \frac{4}{L} \geq 1 - \tfrac{2 \exp(r^2)}{L} + \frac{4}{L} .
    \end{align*}
    When additionally \(\exp(r^2) \leq 3\), then clearly \(A(r) \geq 1 - \frac{2}{L} = 0.8\).
    Thus, it remains to verify the case \(\exp(r^2) > 3\), i.e., \(r > \sqrt{\ln(3)}\).
    We compute the derivative
    \begin{equation*}
        A'(r) = \exp(r^2) \left( 2 - \frac{4r}{L} + 4 r^2 \right) - 4 = \exp(r^2) \left( 2 - \frac{1}{L^2} + \left( \frac{1}{L} -  2 r \right)^2 \right) - 4.
    \end{equation*}
    Since \(\exp(r^2) \geq 3\) this yields
    \begin{align*}
        A'(r) &\geq 3 \left( 2 - \frac{1}{L^2} + \left( \frac{1}{L} -  2 r \right)^2 \right) - 4 = 2 - \frac{3}{L^2} + 3 \left( \frac{1}{L} -  2 r \right)^2 \geq 2 - \frac{3}{100} > 0
    \end{align*}
    meaning that \(A(r)\) is strictly increasing for \(r \geq \sqrt{\ln(3)}\). 
    Since \(A(\sqrt{\ln(3)}) = 1+ 2 \sqrt{\ln(3)} - \frac{2}{L} \geq 0.8\), we conclude that \(A(r) \geq 0.8\) for \(r \geq \sqrt{\ln(3)}\).
    Putting everything together, we have shown that \(A(r) \geq 0.8\) for \(r \geq 0\).

    As for the weight \(B\), remark that \[
        \alpha \Vert \nabla f(x) \Vert = (1 + 2 r \left\vert \exp(r^2) - 2 \right\vert) 2 r \left\vert \exp(r^2) - 2 \right\vert \geq 4 r^2 \left( \exp(r^2) - 2 \right)^2,
    \]
    and \[
        \frac{4 \exp(r^2)}{L} \Vert x \Vert^2 = \frac{4 r^2 \exp(r^2)}{10}.
    \]
    Therefore, the second coefficient is lower bounded by
    \begin{align*}
        B(r) &\geq 4 r^2 \left( \left( \exp(r^2) - 2 \right)^2 - \frac{\exp(r^2)}{10} \right) = 4 r^2 \left( \exp(r^2)^2 - 4 \exp(r^2) + 4 - \frac{\exp(r^2)}{10} \right).
    \end{align*}
    Let us define \(z = r^2 \geq 0\) and \(w = \exp(z) \geq 1\).
    Then we can express this lower bound as
    \begin{equation*}
        Q(w, z) := z \left( 4 w^2 - 16.4 w + 16 \right).
    \end{equation*}
    The quadratic factor is negative only when \(1.6 < w < 2.5 \), and is minimized at \(w^\star = 2.05\) where it attains the minimum value \(-0.81\).
    It remains to lower bound \(Q(w, z)\) for \(\ln(1.6) < z < \ln(2.5)\).
    We have
    \begin{equation*}
        Q(w, z) > -0.81 z > -0.81 \ln(2.5) > -0.8.
    \end{equation*}
    We conclude that \(B(r) \geq -0.8\) for \(r \geq 0\), which completes the proof.
\end{proof}

\subsection{Connection to \texorpdfstring{\((\rho, L_0, L_\rho)\)-smoothness}{(R, L0, LR)-smoothness}}

This section establishes a connection between anisotropic smoothness and \((\rho, L_0, L_\rho)\)-smoothness, which arguably describes the most important subset of \(\ell\)-smooth functions \cite{li_convex_2023}. 
\begin{definition}[{\cite[Definition 3]{li_convex_2023}}]
    A twice continuously differentiable function \(f : \R^n \to \R\) is \((\rho, L_0, L_\rho)\)-smooth for constants \(\rho, L_0, L_\rho \geq 0\) if \(\Vert \nabla^2 f(x) \Vert \leq L_0 + L_\rho \Vert \nabla f(x) \Vert^\rho\) for all \(x \in \R^n\).
\end{definition}
Note that the original definition is slightly more general, as it encompasses functions without full domain and only requires the Hessian upper bound almost everywhere.
\begin{theorem}
    Suppose that a univariate function \(f : \R \to \R\) is \((\rho, L_0, L_\rho)\)-smooth for constants \(\rho, L_0, L_\rho \geq 0\), with \(\rho \leq 2\).
    Then, \(f\) satisfies the second-order characterization of \((L, \bar L)\)-anisotropic smoothness relative to the reference function \(\phi(x) = - \Vert x \Vert - \ln(1 - \Vert x \Vert)\) for \(\bar L = 1\) and \(L \geq 2 \max \{L_0, L_\rho\}\).
\end{theorem}
\begin{proof}
    By \cite[Lemma 2.5]{oikonomidis_nonlinearly_2025}, the second-order characterization of anisotropic smoothness for \(\bar L = 1\) is equivalent to
    \begin{equation*}
        L^{-1} \nabla^2 f(x) \preceq \left[ \nabla^2 \phi^*(\nabla f(x)) \right]^{-1}.
    \end{equation*}
    For our particular reference function and because \(f\) is univariate, this condition becomes
    \begin{equation*}
        L^{-1} f''(x) \leq (1+\vert f'(x) \vert)^2 = 1 + 2 \vert f'(x) \vert + f'(x)^2. 
    \end{equation*}
    We now prove that this upper bound holds.
    By \((\rho, L_0, L_\rho)\)-smoothness we have for all \(x \in \R\) that
    \begin{equation*}
        L^{-1} f''(x) \leq L^{-1} \vert f''(x) \vert \leq \frac{L_0}{L} + \frac{L_\rho}{L} \vert f'(x) \vert^\rho. 
    \end{equation*}
    We distinguish two cases.
    If \(\vert f'(x) \vert \leq 1\), then \(L^{-1} f''(x) \leq \frac{L_0}{L} + \frac{L_\rho}{L}\) and for \(L \geq 2 \max \{L_0, L_\rho\}\) we obtain \(L^{-1} f''(x) \leq 1\), which establishes the required upper bound.
    If \(\vert f'(x) \vert > 1\), then it follows from \(\rho \leq 2\) and \(L \geq \max \{L_0, L_\rho\}\) that
    \begin{equation*}
        L^{-1} f''(x) \leq \frac{L_0}{L} + \frac{L_\rho}{L} \vert f'(x) \vert^2 \leq 1 + \vert f'(x) \vert^2,
    \end{equation*}
    which implies the required upper bound.
\end{proof}
Remark that most examples of univariate \((\rho, L_0, L_\rho)\)-smooth functions in \cite{li_convex_2023} satisfy \(\rho \leq 2\).
Besides polynomials, this includes exponential functions \(a^x\) with \(a > 1\), and double exponentials \(a^{(b^x)}\) with \(a, b > 0\).

\subsection{Verification of assumptions on the reference functions} \label{sec:assumptions}

Throughout this work, we have made a number of assumptions which only relate to the reference function \(\phi\), i.e., \cref{assump:basic-reference,assumption:phi-star-c2,assumption:isotropic,assump:perturbed}.
This section explicitly verifies these assumptions for \emph{isotropic} reference functions \(\phi = h \circ \Vert \cdot \Vert\) where the kernel function \(h\) is one of the following:
\begin{equation*}
    h_1(x) = \cosh(x) - 1, \quad h_2(x) = \exp(\vert x \vert) - \vert x \vert - 1, \quad h_3(x) = - \vert x \vert - \ln(1 - \vert x \vert).
\end{equation*}
To this end, the following results from \cite[Table 1]{oikonomidis_nonlinearly_2025} will prove useful:
\begin{align*}
    &{h_1^*}'(y) = \arcsinh(y), \quad &&{h_2^*}'(y) = \ln(1+\vert y \vert) \barsgn(y), \quad &&{h_3^*}'(y) = \tfrac{y}{1+\vert y \vert}\\
    &{h_1^*}''(y) = \tfrac{1}{\sqrt{1+y^2}}, \quad &&{h_2^*}''(y) = \tfrac{1}{1+\vert y \vert}, \quad &&{h_3^*}''(y) = \tfrac{1}{(1+\vert y \vert)^2}.
\end{align*}

\paragraph{Assumptions \ref{assump:basic-reference} and \ref{assumption:phi-star-c2}}

These assumptions were proven for \(h_1, h_2, h_3\) in \cite{oikonomidis_nonlinearly_2025}.

\paragraph{Assumption \ref{assumption:isotropic}}

To verify (i), i.e., whether \(\tfrac{{h^*}'(y)}{y} \) is decreasing on \(\R_+\), we can check if for \(y > 0\)
\begin{equation*}
    \tfrac{d}{dy} \left( \tfrac{{h^*}'(y)}{y} \right) = \tfrac{y {h^*}''(y) - {h^*}'(y)}{y^2} < 0,.
\end{equation*}
or equivalently, \(y {h^*}''(y) < {h^*}'(y)\).
This holds for \(h_1\), \(h_2\), and \(h_3\). 
Part (ii) holds, since
\begin{equation*}
    \lim_{y \to \infty} y {h_1^*}''(y) = 1, \qquad \lim_{y \to \infty} y {h_2^*}''(y) = 1, \qquad \lim_{y \to \infty} y {h_3^*}''(y) = 0.
\end{equation*}
Also (iii) is satisfied, since \( {h_1^*}'\) and \( {h_2^*}'\) scale logarithmically and \( {h_3^*}'\) is bounded.
In particular,
\begin{align*}
    \lim_{y \to +\infty} \tfrac{{h_1^*}'(s_d(y))}{y} &= \lim_{y \to +\infty} \tfrac{d \ln(y)}{y} = 0\\
    \lim_{y \to +\infty} \tfrac{{h_2^*}'(s_d(y))}{y} &= \lim_{y \to +\infty} \tfrac{d \ln(y)}{y} = 0\\
    \lim_{y \to +\infty} \tfrac{{h_3^*}'(s_d(y))}{y} &= \lim_{y \to +\infty} \tfrac{1}{y} = 0.
\end{align*}
Here we used the fact that \(\arcsinh(y) = \ln(y + \sqrt{y^2 + 1})\).
Thus, \cref{assumption:isotropic} holds for \(h_1\), \(h_2\) and \(h_3\).

\paragraph{Assumption \ref{assump:perturbed}}

The kernel functions have the following Taylor expansions:
\begin{align*}
    h_1(x) &= \sum_{n = 1}^\infty \tfrac{x^{2n}}{(2n)!} = \tfrac{x^2}{2!} + \tfrac{x^4}{4!} + \dots\\
    h_2(x) &= \sum_{n = 2}^\infty \tfrac{\vert x \vert^{n}}{n!} = \tfrac{\vert x \vert^2}{2!} + \tfrac{|x|^3}{3!} + \dots\\
    h_3(x) &= \sum_{n = 2}^\infty \tfrac{\vert x \vert^{n}}{n} = \tfrac{\vert x \vert^2}{2} + \tfrac{|x|^3}{3} + \dots
\end{align*}
We remark in particular that each summand is nonnegative, and that the term with the lowest degree equals \(\nicefrac{x^2}{2}\) for all three kernel functions.
This immediately implies \(h_i(x) \geq \nicefrac{x^2}{2}\) and \(h_i(x) = \nicefrac{x^2}{2} + o(x^2)\) as \(x \to 0\), for \(i \in \{1, 2, 3\}\), and proves that the assumption holds for \(h_1\), \(h_2\) and \(h_3\).

\subsection{Generalized smoothness of regularized neural networks with quadratic loss}

The main goal of this section is to investigate the \emph{generalized smoothness} of the quadratic loss of a deep neural network, and in particular, whether \cref{assump:novel-condition} holds.
It turns out that, under sufficient regularization, this is indeed the case.
To establish this result, we need to bound the gradient and Hessian norm of the loss by a polynomial of an appropriate degree.

Consider a deep \(N\)-layer neural network with quadratic loss.
Each layer -- of which we denote the index by \(t \in \N_{[0, N-1]}\) -- consists of weights \(W_t \in \R^{n_{t+1} \times n_{t}}\), a bias term \(b^t \in \R^{n_{t+1}}\), and a componentwise activation function \(\Sigma_t : \R^{n_{t+1}} \to \R^{n_{t+1}}\).
The final mapping \(\Sigma_{N-1}\) is assumed to be the identity mapping, as common in regression problems.
For a given data point \(\bar x \in \R^{n_0}\) and corresponding label \(\bar y \in \R^{n_{N}}\), training this network entails minimizing the loss function
\begin{equation} \label{eq:nn-loss}
    \frac{1}{2} \Vert \Sigma_{N-1}(W_{N-1} \Sigma_{N-2}(\dots \Sigma_0(W_0 \bar x + b^0) \dots ) + b^t) - \bar y \Vert^2.
\end{equation}
Neural networks are usually trained on a large set of pairs \(\{\bar x^i, \bar y^i \}_{i \in \N_{[1,I]}}\), in which case the total loss becomes a summation of the losses for each individual pair.
To simplify the presentation, we proceed with \(I = 1\), i.e., with \eqref{eq:nn-loss}, but the results are easily extended to the case \(I > 1\). 

Let us denote intermediate variables \(x^{t+1} = f_t(x^t, w^t) = \Sigma_t(W_t x^t + b^t) \in \R^{n_{t+1}}\), where we use the convention \(x^0 = \bar x\), and define the vectorized weights and bias at layer \(t\) by \(w^t = \left(
    \vecop(W_t)^\top, {b^t}^\top
\right)^\top \in \R^{(n_t+1)n_{t+1}}\).
Then the loss \eqref{eq:nn-loss} can be compactly written as \(\frac{1}{2} \Vert x^N - y \Vert^2\).

The minimization of the loss \eqref{eq:nn-loss} can then be interpreted as an optimal control problem (OCP) of horizon \(N\) with states \(x^t \in \R^{n_t}\) and inputs \(w^t \in \R^{(n_t+1)n_{t+1}}\), i.e.,
\begin{equation*}
    \begin{aligned}
        \minimize{}\quad &\sum_{t = 0}^{N-1} \ell_t(x^t, w^t) + \ell_N(x^N) && \stt{}\quad &&x^{t+1} = f_t(x^t, w^t), t = 0, \dots, N-1,
    \end{aligned}
\end{equation*}
where \(x^0 = \bar x\), \(\ell_t \equiv 0\) for \(t \in \N_{[0,N-1]}\) and
\(
    \ell_N(x^N) = \frac{1}{2} \Vert x^N - \bar y \Vert^2
\).
In the context of OCPs, the functions \(f_t(x^t, w^t) = \Sigma_t(W_t x^t + b^t)\) are called the dynamics, and have gradients
\begin{equation} \label{eq:nn-dynamics-gradient}
    \begin{aligned}
        \nabla_{x^t} f_t(x^t, w^t) &= W_t^\top (\Sigma_t'(W_t x^t + b^t))^\top &&\in \R^{n_t \times n_{t+1}}\\
        \nabla_{w^t} f_t(x^t, w^t) &= \begin{pmatrix}
            x^t \otimes I_{n_{t+1}}\\
            I_{n_{t+1}}
        \end{pmatrix}(\Sigma_t'(W_t x^t + b^t))^\top &&\in \R^{(n_t + 1) n_{t+1} \times n_{t+1}}.
    \end{aligned}
\end{equation}
Let \(x := ({x^1}^\top, \dots, {x^N}^\top)^\top\) and \(w := ({w^0}^\top, \dots, {w^{N-1}}^\top)^\top\) denote the vectors containing all states and inputs respectively.
We aim to derive expressions for the gradient and Hessian of the loss function \eqref{eq:nn-loss} with respect to \(w\).
To that end, we use a standard idea in optimal control and eliminate the dynamics. 
This approach is known as \emph{single shooting}.
Following Bertsekas~\cite[\S 1.9]{bertsekas_nonlinear_1999}, we introduce mappings
\begin{equation*}
    F_{t+1}(w) = f_t(F_t(w), w) = x^{t+1}, \quad t = 0, \dots, N-1.
\end{equation*}
with \(F_{-1}(w) = \bar x\).
Let \(F(w) = (F_1(w)^\top, \dots, F_N^\top)^\top\) and note that \cite[Eq. (1.246) and below]{bertsekas_nonlinear_1999}
\begin{equation} \label{eq:nabla-F}
    \nabla F(w) = \begin{pmatrix}
        \nabla_{w^0} f_0 & \nabla_{w^0} f_0 \nabla_{x^1} f_1  & \dots & \nabla_{w^0} f_0 \nabla_{x^1} f_1 \dots \nabla_{x^{N-1}} f_{N-1}\\
        0 & \nabla_{w^1} f_1 & \dots & \nabla_{w^1} f_1 \nabla_{x^2} f_2 \dots \nabla_{x^{N-1}} f_{N-1}\\
        \vdots & \vdots & \dots & \vdots\\
        0 & 0 & \dots & \nabla_{w^{N-1}} f_{N-1}
    \end{pmatrix} %
\end{equation}
We highlight that the block columns correspond to the gradients of the individual mappings \(F_t\), i.e.,\begin{equation*}
    \nabla F(w) = \begin{pmatrix}
        \nabla F_1(w) & \nabla F_2(w) & \dots & F_N(w)
    \end{pmatrix}.
\end{equation*}
Therefore, the OCP, and equivalently, the neural network training problem, is compactly written as
\begin{equation} \label{eq:nn-loss-J}
    \minimize_{w} J(w) := \ell_N(F_N(w)) = \frac{1}{2} \Vert F_N(w) - \bar y \Vert^2
\end{equation} 
The gradient of \(J\) is then easily expressed by 
\begin{equation} \label{eq:nn-loss-grad}
    \nabla J(w) = \nabla_w F_N(w) (F_N(w) - y).
\end{equation}
Obtaining an expression for the Hessian is more involved.
We follow \cite[\S 1.9]{bertsekas_nonlinear_1999}, which uses the Lagrangian function of the OCP with multipliers \(\lambda = ({\lambda^1}^\top, \dots, {\lambda^N}^\top)^\top \in \R^{\sum_{t = 0}^{N-1} n_{t+1}}\), i.e.,
\begin{equation*}
    L(x, w, \lambda) = \ell_N(x) + \sum_{t = 0}^{N-1} (f_t(x^t, w^t) - x^{t+1})^\top \lambda^{t+1}.
\end{equation*}
The central idea is to express \(\nabla J(w)\) in terms of \(\nabla_x L(x, w, \lambda)\) and \(\nabla_w L(x, w, \lambda)\).
The computation of \(\nabla^2 J(w)\) is simplified by selecting an appropriate multiplier, which is recursively defined by
\begin{equation*}
    \begin{aligned}
        \lambda^N &= F_N - \bar y && \in \R^{n_N}\\
        \lambda^t &= (\nabla_{x^t} f_t) \lambda^{t+1} &&\in \R^{n_{t}}, \qquad t = 1, \dots N-1.
    \end{aligned}
\end{equation*} 
This yields for example \(\lambda^1 = \nabla_{x^1} f_1 \nabla_{x^2} f_2 \dots \nabla_{x^{N-1}} f_{N-1} (F_N - y)\).
Bertsekas~\cite[Eq. (1.242)]{bertsekas_nonlinear_1999} establishes that
\begin{equation} \label{eq:nn-loss-hessian}
    \begin{aligned}
        \nabla^2 J(w) = \nabla F(w) \nabla_{xx}^2 L(F(w), w, \lambda) \nabla F(w)^\top &+ 2 \nabla F(w) \nabla_{xw}^2 L(F(w), w, \lambda)\\
        &+ \nabla_{ww}^2 L(F(w), w, \lambda).
    \end{aligned}
\end{equation}

The following theorem establishes that, under sufficient regularization, the regularized loss function of neural network training satisfies \cref{assump:novel-condition}, and consequently satisfies the generalized smoothness notions investigated in this work.
We emphasize that the rather restrictive bound on the power \(P\) of the regularizer can be relaxed significantly by exploiting the structure of \(\nabla F(w)\) (cf.\,\eqref{eq:nabla-F}) and further working out the Hessian \(\nabla^2 J(w)\) (cf.\,\eqref{eq:nn-loss-hessian}) before upper bounding its norm.
However, for simplicity, the proof below immediately uses submultiplicativity of the matrix norm (cf.\,\eqref{eq:nn-hessian-norm-bound}).  
\begin{theorem}
    Consider the regularized neural network training problem with objective 
    \begin{equation*}
        \tilde J(w) = J(w) + \kappa \Vert w \Vert^{P}
    \end{equation*}
    where \(J(w)\) is the quadratic loss of an \(N\)-layer network as defined in \eqref{eq:nn-loss-J}.
    Suppose that the mappings \(\Sigma_t\) are bounded and have bounded first and second derivatives for \(t \in \N_{[0, N-2]}\), and that the mapping \(\Sigma_{N-1}\) is the identity map.
    If \(P \geq 3N+2\), then the following statements hold for any \(\kappa > 0\).
    \begin{theoremenum}
        \item For any \( L_1 > 0\) there exists an \(L_0 > 0\) such that \(\tilde J\) is \((L_0, L_1)\)-smooth.
        \item If \(\phi\) satisfies \cref{assumption:phi-star-c2,assumption:isotropic}, then for any \(\bar L > 0\), there exists an \(L > 0\) such that \(\tilde J\) satisfies the second-order characterization of \((L, \bar L)\)-anisotropic smoothness.
    \end{theoremenum}
\end{theorem}
\begin{proof}
By boundedness of \(\Sigma_{t}\), also the states \(x^{t+1}\) are bounded for \(t \in \N_{[0,N-2]}\), whereas we note that \(\Vert x^N \Vert = \Vert F_N(w)  \Vert = \cO(\Vert w \Vert)\).
From the gradient expressions \eqref{eq:nn-dynamics-gradient}, we observe that \(\Vert \nabla_{w^t} f_t \Vert = \cO(1)\) and \(\Vert \nabla_{x^t} f_t \Vert = \cO (\Vert w \Vert)\), and by \eqref{eq:nabla-F} we obtain \(\Vert \nabla F(w) \Vert = \cO(\Vert w \Vert^{N-1})\).
Therefore, the gradient \(\nabla J(w)\) as defined in \eqref{eq:nn-loss-grad} is upper bounded by a polynomial of degree \(N\) in \(\Vert w \Vert\), i.e.,
\(
    \Vert \nabla J(w) \Vert = \cO (\Vert w \Vert^{N}).
\)
Consequently,
\begin{equation} \label{eq:nn-reg-loss-gradient}
    \Vert \nabla \tilde J(w) \Vert \geq \kappa P \Vert w \Vert^{P-1} - \Vert \nabla J(w) \Vert = \kappa P \Vert w \Vert^{P-1} - \cO(\Vert w \Vert^{N}).
\end{equation}
As for the Hessian norm, we have by \eqref{eq:nn-loss-hessian} that
\begin{equation} \label{eq:nn-hessian-norm-bound}
    \Vert \nabla^2 J(w) \Vert = \cO \left( \Vert w \Vert^{2N-2} \bigg \Vert \begin{pmatrix}
        \nabla_{ww}^2 L(F(w), w, \lambda) & \nabla_{wx}^2 L(F(w), w, \lambda)\\
        \nabla_{xw}^2 L(F(w), w, \lambda) & \nabla_{xx}^2 L(F(w), w, \lambda)
    \end{pmatrix} \bigg\Vert \right).
\end{equation}
We compute the gradient of the Lagrangian with respect to \(x\) and \(w\)
\begin{align*}
    \nabla_{x^t} L(x, w, \lambda) &= (\nabla_{x^t} f_t(x^t, w^t)) \lambda^{t+1} - \lambda^t \qquad &&t = 1, \dots N-1\\
    \nabla_{x^N} L(x, w, \lambda) &= (x^N - \bar y) - \lambda^N\\
    \nabla_{w^t} L(x, w, \lambda) &= (\nabla_{w^t} f_t(x^t, w^t)) \lambda^{t+1}, \qquad &&t = 0, \dots, N-1.
\end{align*}
Clearly, \(\nabla_{xx}^2 L(F(w), w, \lambda)\) is block diagonal.
The last block
\(
    \nabla_{x^N x^N}^2 L(x, w, \lambda) = I_{n_N}
\)
is straightforward to compute. For the other ones, we proceed by rewriting
\begin{align*}
    \nabla_{x^t} L(x, w, \lambda) &= \vecop((\nabla_{x^t} f_t(x^t, w^t)) \lambda^{t+1}) - \lambda^t = \left( {\lambda^{t+1}}^\top \otimes I_{n_{t}} \right) \vecop((\nabla_{x^t} f_t(x^t, w^t))) - \lambda^t.
\end{align*}
Here we used the fact that \(\vecop(z) = z\) for any vector, and \(
    \vecop(I B A) = (A^\top \otimes I) \vecop(B)
\)
for any matrices \(A, B\) and identity matrix \(I\) of compatible dimensions.
This yields
\begin{equation*}
    \nabla_{x^t x^t}^2 L(x, w, \lambda) = \left( {\lambda^{t+1}}^\top \otimes I_{n_{t}} \right) \frac{\partial \vecop((\nabla_{x^t} f_t(x^t, w^t)))}{\partial x^t}.
\end{equation*}
From the fact that \(\Vert \lambda^{t+1} \Vert = \cO(\Vert w \Vert^{N-t})\) for \(t = 0, \dots, N-1\), and from boundedness of \(\Sigma_t''\) we obtain 
\begin{equation*}
    \Vert \nabla_{x^t x^t}^2 L(x, w, \lambda) \Vert = \cO (\Vert w \Vert^{N-t+2}), \qquad t = 0, \dots, N-1.
\end{equation*}
and it follows that 
\(
    \Vert \nabla_{xx}^2 L(x, w, \lambda) \Vert = \cO(\Vert w \Vert^{N+2}).
\)
Similar arguments can be used to show that also the other Hessian blocks of the Lagrangian satisfy a similar upper bound.
In conclusion, we obtain that 
\(
    \Vert \nabla^2 J(w) \Vert = \cO (\Vert w \Vert^{2N-2} \Vert w \Vert^{N+2}) = \cO (\Vert w \Vert^{3N}).
\)
It follows immediately that \(\Vert \nabla^2 \tilde J (w) \Vert\) is upper bounded by a polynomial of degree \(3N\) in \(\Vert w \Vert\).
And since by \eqref{eq:nn-reg-loss-gradient} the gradient \(\Vert \nabla \tilde J(w) \Vert\) is lower bounded by a polynomial of degree \(3N+1\) in \(\Vert w \Vert\) with strictly positive leading coefficient, we conclude that \cref{assump:novel-condition} holds.
The claims then follow by \cref{thm:l0-l1-smoothness-polynomials,thm:anisotropic-smoothness-polynomials}.
\end{proof}

\section{Auxiliary results}
\begin{lemma}[Power Mean inequality] \label{lem:power-mean}
Let \( p > q > 0 \). Then,
    \begin{equation*}
        \left( \frac{1}{m} \sum_{i = 1}^m \vert \alpha_i \vert^p \right)^{\nicefrac{1}{p}} \geq \left( \frac{1}{m} \sum_{i = 1}^m \vert \alpha_i \vert^q\right)^{\nicefrac{1}{q}}.
    \end{equation*}
\end{lemma}
\begin{proof}
    Note that the function \( \varphi(\alpha) = \alpha^{\nicefrac{q}{p}} \) is concave for \(\alpha > 0\) since \(q < p\).
    By Jensen's inequality, this implies
    \begin{equation*}
        \left( \frac{1}{m} \sum_{i = 1}^m \vert \alpha_i \vert^p \right)^{\nicefrac{q}{p}} = \varphi \left( \frac{1}{m} \sum_{i = 1}^m \vert \alpha_i \vert^p \right) \geq \frac{1}{m} \sum_{i = 1}^m \varphi(\vert \alpha_i \vert^p) = \frac{1}{m} \sum_{i = 1}^m \vert \alpha_i \vert^q.
    \end{equation*}
    Raising both sides to the power \(\nicefrac{1}{q}\) establishes the claim.
\end{proof}

\begin{lemma} \label{lem:squared-norm-to-norm-lower-bound}
Let \( F : \mathbb{R}^n \to \mathbb{R}^m \) be a mapping satisfying, for all \( x \in \mathbb{R}^n \),
\[
\|F(x)\|^2 \geq a \|x\|^6 - b \|x\|^4,
\]
where \( a > 0 \) and \( b \geq 0 \). Then, for all \( x \in \mathbb{R}^n \),
\[
\|F(x)\| \geq \sqrt{a} \|x\|^3 - \frac{b}{\sqrt{a}} \|x\|.
\]
\end{lemma}

\begin{proof}  
In the proof we write \( r = \|x\| \) for simplicity. The given inequality becomes:
\[
\|F(x)\|^2 \geq a r^6 - b r^4 = r^4(a r^2 - b).
\]
We distinguish two cases:

\textit{Case 1:} \( r^2 \geq \frac{b}{a} \).  
In this case the given lower bound is nonnegative. Taking the square root yields
\[
\|F(x)\| \geq \sqrt{r^4 (a r^2 - b)} = r^2 \sqrt{a r^2 - b}.
\]
We apply the inequality \( \sqrt{1 - z} \geq 1 - z \) for \( 0 \leq z \leq 1 \), by letting \( z = \frac{b}{a r^2} \in [0,1] \). Then
\[
\sqrt{a r^2 - b}
= \sqrt{a r^2 \left(1 - \frac{b}{a r^2} \right)}
= \sqrt{a} r \sqrt{1 - \frac{b}{a r^2}}
\geq \sqrt{a} r \left(1 - \frac{b}{a r^2} \right)
= \sqrt{a} r - \frac{b}{\sqrt{a} r}.
\]
Multiplying both sides by \( r^2 \), we obtain
\[
\|F(x)\| \geq r^2 \left( \sqrt{a} r - \frac{b}{\sqrt{a} r} \right)
= \sqrt{a} r^3 - \frac{b}{\sqrt{a}} r.
\]

\textit{Case 2:} \( r^2 < \frac{b}{a} \).
The claim follows immediately by nonnegativity of the norm and from the fact that the desired lower bound is negative in this case, i.e.,
\begin{equation*}
    \Vert F(x) \Vert \geq 0 > \frac{r}{\sqrt{a}} \left( a r^2 - b \right) = \sqrt{a} r^3 - \frac{b}{\sqrt{a}} r.
\end{equation*}
\end{proof}

\begin{lemma} \label{lem:upper-bound-preconditioned-gradient-norm}
    Suppose that \cref{assump:perturbed} holds, and that a point \(x \in \R^n\) satisfies
    \begin{equation*}
        \lambda^{-1} \phi(\nabla \phi^*(\lambda \nabla f(x))) \leq \frac{\cG^2}{2},
    \end{equation*}
    where \(\cG\) is defined as in \eqref{eq:perturbed-GD-parameters-1}.
    Then,
    \begin{equation*}
        \Vert \nabla \phi^*(\lambda \nabla f(x)) \Vert \leq \sqrt{\lambda} \cG.
    \end{equation*}
\end{lemma}
\begin{proof}
    The bound \(\phi(x) \geq \frac{\Vert x \Vert^2}{2}\) yields 
    \begin{equation*}
        \begin{aligned}
            \phi(\nabla \phi^*(\lambda \nabla f(x))) \leq \frac{\lambda \cG^2}{2} \leq \phi(\sqrt{\lambda} \cG).
        \end{aligned}
    \end{equation*}
    The claim then follows by \(\phi = h \circ \Vert \cdot \Vert\), nonnegativity of \(\cG\) and strict monotonicity of \(h\).
\end{proof}

\section{Missing proofs of section \ref{sec:anisotropic-smoothness}}
\subsection{Proof theorem \ref{prop:second-order-characterization}}

\begin{proof}
    This proposition is a direct combination of \cite[Propositions 2.6 \& 2.9]{oikonomidis_nonlinearly_2025}.
\end{proof}

\subsection{Proof of theorem \ref{thm:l0-l1-smoothness-polynomials}}

\begin{proof}    
    By \cref{assump:novel-condition} we have
    \begin{equation*}
        \limsup_{\Vert x \Vert \to \infty} \frac{\Vert \nabla^2 f(x) \Vert_F}{\Vert \nabla f(x) \Vert} \leq \limsup_{\Vert x \Vert \to \infty} \frac{p_{R}(\Vert x \Vert)}{q_{R+1}(\Vert x \Vert)} = \limsup_{\Vert x \Vert \to \infty} \frac{a_R}{b_{R+1} \Vert x \Vert} = 0.
    \end{equation*}
    By nonnegativity of \(\nicefrac{\Vert \nabla^2 f(x) \Vert_F}{\Vert \nabla f(x) \Vert}\) we conclude that \(\lim_{\Vert x \Vert \to \infty} \frac{\Vert \nabla^2 f(x) \Vert_F}{\Vert \nabla f(x) \Vert} = 0\). Thus, for any \( L_1 > 0\) there exists \( \delta > 0\) such that
    \begin{equation*}
        \Vert x \Vert > \delta \Rightarrow \Vert \nabla^2 f(x) \Vert_F \leq L_1 \Vert \nabla f(x) \Vert.
    \end{equation*}
    Moreover, by continuity of \(\nabla^2 f\), we know that \( \Vert \nabla^2 f(x) \Vert_F \) is bounded on the compact set \(\Omega := \{ x \mid \Vert x \Vert \leq \delta \}\).
    We conclude that \(f\) is \((L_0, L_1)\)-smooth with \(L_0 = \max_{x \in \Omega} \Vert \nabla^2 f(x) \Vert_F\).
\end{proof}

\subsection{Proof of theorem \ref{thm:anisotropic-smoothness-polynomials}}

\begin{proof}
    Fix an arbitrary \( \bar L > 0\) and, for ease of notation, define
    \begin{equation*}
        H_\lambda(x) := \nabla^2 \phi^*(\bar L^{-1} \nabla f(x)) \nabla^2 f(x).
    \end{equation*}
    If \(\lim_{\Vert x \Vert \to \infty} \Vert H_\lambda(x) \Vert = 0\), then for any \(\epsilon > 0 \) there exists \(\delta > 0\) such that
    \begin{equation*}
        \Vert x \Vert > \delta \Rightarrow \Vert H_\lambda(x) \Vert < \epsilon.
    \end{equation*}
    Therefore, the continuous function \(x \to \Vert H_\lambda(x) \Vert\) is bounded on the compact set \(\Omega := \{ x \mid \Vert x \Vert \leq \delta \}\), and for all \( x \in \R^n \setminus \Omega\) we know that \(\Vert H_\lambda(x) \Vert \leq \epsilon\).
    We conclude that if \(\lim_{\Vert x \Vert \to \infty} \Vert H_\lambda(x) \Vert = 0\), then \(\Vert H_\lambda(x) \Vert\) is bounded on \(\R^n\), and because \(\lambda_{\max}(H_\lambda(x)) \leq \Vert H_\lambda(x) \Vert \), this would prove the claim.
    By equivalence of norms, it suffices to show boundedness of any norm; we proceed with the Frobenius norm \(\Vert H_\lambda(x) \Vert_F \).
    Because for isotropic reference functions
    \begin{equation*}
        \nabla^2 \phi^*(y) = {h^*}''(\Vert y \Vert) \frac{y y^\top}{\Vert y \Vert^2} + \frac{{h^*}'(\Vert y \Vert)}{\Vert y \Vert} \left( I - \frac{y y^\top}{\Vert y \Vert^2} \right),
    \end{equation*}
    it follows that
    \begin{equation*}
        \begin{aligned}
            \Vert H_\lambda(x) \Vert_F &\leq \bigg \Vert {h^*}''({\scriptstyle \Vert \bar L^{-1} \nabla f(x)\Vert}) \tfrac{\nabla f(x) \nabla f(x)^\top}{\Vert \nabla f(x) \Vert^2} + \tfrac{{h^*}'(\Vert \bar L^{-1} \nabla f(x)\Vert )}{\Vert \bar L^{-1} \nabla f(x) \Vert} \left( I - \tfrac{\nabla f(x) \nabla f(x)^\top}{\Vert \nabla f(x) \Vert^2} \right) \bigg \Vert \Vert \nabla^2 f(x) \Vert_F\\
            &\leq \left( \vert {h^*}''({\scriptstyle \Vert \bar L^{-1} \nabla f(x) \Vert}) \vert + \tfrac{\vert {h^*}'(\Vert \bar L^{-1} \nabla f(x) \Vert) \vert}{\Vert \bar L^{-1} \nabla f(x) \Vert} \right) p_R(\Vert x \Vert).
        \end{aligned}
    \end{equation*}
    Taking the limes superior and using \cref{assumption:isotropic} yields
    \begin{equation*}
        \begin{aligned}
            \limsup_{\Vert x \Vert \to \infty} \Vert H_\lambda(x) \Vert_F &= \limsup_{\Vert x \Vert \to \infty} \left( \vert {h^*}''(\Vert \bar L^{-1} \nabla f(x) \Vert) \vert + \frac{\vert {h^*}'(\Vert \bar L^{-1} \nabla f(x) \Vert) \vert}{\Vert \bar L^{-1} \nabla f(x) \Vert} \right) p_R(\Vert x \Vert)\\
            &=\limsup_{\Vert x \Vert \to \infty} \left( \frac{C_2}{\Vert \bar L^{-1} \nabla f(x) \Vert} + \frac{\vert {h^*}'(\Vert \bar L^{-1} \nabla f(x) \Vert) \vert}{\Vert \bar L^{-1} \nabla f(x) \Vert} \right) p_R(\Vert x \Vert)\\
            &\leq \limsup_{\Vert x \Vert \to \infty} \left( \frac{C_2}{\bar L^{-1} q_{R+1}(\Vert x \Vert)} + \frac{\vert {h^*}'\left(\bar L^{-1} q_{R+1}(\Vert x \Vert)\right) \vert}{\bar L^{-1} q_{R+1}(\Vert x \Vert)} \right) p_R(\Vert x \Vert)\\
            &\leq \limsup_{\Vert x \Vert \to \infty} \frac{a_R C_2}{b_{R+1} \Vert x \Vert} + \frac{a_R \vert {h^*}'\left(\bar L^{-1} q_{R+1}(\Vert x \Vert)\right) \vert}{b_{R+1} \Vert x \Vert} = 0.
        \end{aligned}
    \end{equation*}
    Here, the third step used the fact that \(\nicefrac{{h^*}'(x)}{x}\) is decreasing on \( \R_+\) in combination with \(\Vert \bar L^{-1} \nabla f(x) \Vert \geq \bar L^{-1} q_{R+1}(\Vert x \Vert)\).
    Since \( \Vert \cdot \Vert \geq 0 \), we conclude \( \lim_{\Vert x \Vert \to \infty} \Vert H_\lambda(x) \Vert = 0\).
\end{proof}

\subsection{Proof of theorem \ref{th:phase-retrieval-smoothness}}

\begin{lemma}[Gradient norm lower bound] \label{lemma:phase-retrieval-gradient-norm}
    Consider the phase retrieval problem with objective \eqref{eq:phase-retrieval-objective-gradient} and suppose that that the vectors \(\{a_i\}_{i = 1}^m\) span \(\R^n\).
    Then, there exists a constant \(C > 0\) which depends on the measurement vectors \(\{a_i\}_{i = 1}^m\) (but not on \(x\)) such that
    \begin{equation*}
        \Vert \nabla f(x) \Vert \geq C \Vert x \Vert^3 - \sum_{j = 1}^m \Vert y_j^2 a_j a_j^\top \Vert \Vert x \Vert
    \end{equation*}
\end{lemma}
\begin{proof}
    Clearly, the gradient norm can be lower bounded by
    \begin{equation*}
    \begin{aligned}
        \Vert \nabla f(x) \Vert &\geq \bigg \Vert \sum_{i = 1}^m (a_i^\top x)^3 a_i \bigg \Vert - \sum_{j = 1}^m \Vert y_j^2 a_j a_j^\top \Vert \Vert x \Vert.
    \end{aligned}
    \end{equation*}
    Let us further lower bound \(\Vert \sum_{i = 1}^m (a_i^\top x)^3 a_i \Vert\) in terms of \( \Vert x \Vert \).
    For ease of notation, we denote \(g(x) := \Vert G(x) \Vert\) where
    \(
        G(x) := \sum_{i = 1}^m (a_i^\top x)^3 a_i.
    \)
    Observe that \(g\) is positively homogeneous of degree \(3\) \cite[Definition 13.4]{rockafellar_variational_1998}, since for any \(\lambda > 0\)
    \begin{equation*}
        g(\lambda x) = \bigg \Vert \sum_{i = 1}^m (\lambda a_i^\top x)^3 a_i \bigg \Vert = \lambda^3 \bigg \Vert \sum_{i = 1}^m (a_i^\top x)^3 a_i \bigg \Vert = \lambda^3 g(x).
    \end{equation*}
    For any \(x \in \R^n\) we have \(
        g(x) = \Vert x \Vert^3 g\left(\nicefrac{x}{\Vert x \Vert}\right),
    \)
    and hence it suffices to lower bound \(g\) on the unit sphere \( \bS^{n-1} := \left\{x \in \R^n \mid \Vert x \Vert = 1 \right\}\).
    Remark that \(g\) achieves a minimum over \(\bS^{n-1}\) because \(\bS^{n-1}\) is compact and \(g\) is continuous.
    Denote this minimum by \(
        C := \min_{u \in \bS^{n-1}} g(u),
    \)
    such that \[
        g(x) \geq C \Vert x \Vert^3.
    \]
    If \( C > 0\), then the proof is done.
    For this reason, we show by contradiction that the case \( C = 0\) is not possible.
    If \(C = 0\), there exists \(u^\star \in \bS^{n-1}\) such that 
    \begin{equation*}
        G(u^\star) = \sum_{i = 1}^m (a_i^\top u^\star)^3 a_i = 0.
    \end{equation*}
    It follows that
    \begin{equation*}
        {u^\star}^\top G(u^\star) = \sum_{i = 1}^m (a_i^\top u^\star)^4 = 0,
    \end{equation*}
    and since \((a_i^\top u^\star)^4 \geq 0\) we conclude that \(a_i^\top u^\star = 0\) for \(i = 1, \dots m\), i.e., \(u^\star\) is orthogonal to all measurement vectors.
    Therefore, \(u^\star \in (\laspan \left\{ a_1, \dots, a_m \right\})^\top = (\R^n)^\top = \{ 0 \}\).
    This results in a contradiction, since \(u^\star = 0 \notin \bS^{n-1}\).
    We conclude that \(C > 0\), which proves the claim.
\end{proof}

\begin{proof}
    The objective \(f\) is a fourth-order polynomial in the variables \(x\).
    Hence, one easily verifies that 
    \begin{equation*}
        \Vert \nabla^2 f(x) \Vert \leq \left(3 \sum_{i = 1}^m \Vert a_i \Vert^4 \right) \Vert x \Vert^2 + \left(\sum_{j = 1}^m y_j^2 \Vert a_j \Vert^2 \right),
    \end{equation*}
    where the upper bound is a polynomial of degree \(2\) in \(\Vert x \Vert\).
    Moreover, by \cref{lemma:phase-retrieval-gradient-norm} the gradient norm \(\Vert \nabla f(x) \Vert \) can be lower bounded by a polynomial of degree \(3\) in \(\Vert x \Vert\), with strictly positive leading coefficient.
    Therefore, \cref{assump:novel-condition} holds.
    \Cref{th:phase-retrieval-smoothness-1} now follows directly from \cref{thm:l0-l1-smoothness-polynomials}, and \cref{th:phase-retrieval-smoothness-2} from \cref{thm:anisotropic-smoothness-polynomials}.
\end{proof}

\subsection{Proof of theorem \ref{th:matrix-factorization-symmetric-smoothness}}

\begin{lemma}[Gradient norm lower bound] \label{lem:matrix-factorization-symmetric-lower-bound}
    Consider the matrix factorization problem with objective \eqref{eq:symmetric-factorization-objective-gradient}.
    Then the gradient norm can be lower bounded by
    \begin{equation*}
        \Vert \nabla f(U) \Vert_F^2 \geq \frac{1}{n^2} \Vert U \Vert_F^6 - 2 \Vert Y \Vert_F \Vert U \Vert_F^4.
    \end{equation*}
\end{lemma}
\begin{proof}
    We have that
    \begin{equation*}
    \begin{aligned}
        \Vert \nabla f(U) \Vert_F^2 &= \trace(\nabla f(U)^\top \nabla f(U)) = \trace(U^\top (U U^\top - Y)(U U^\top - Y)U)\\
        &= \trace(U^\top U U^\top U U^\top U - 2 U^\top Y U U^\top U + U^\top Y Y U)\\
        &= \trace(U^\top U U^\top U U^\top U) - 2 \trace(U^\top Y U U^\top U) + \trace(U^\top Y Y U)\\
        &\geq \trace(U^\top U U^\top U U^\top U) - 2 \Vert Y \Vert_F \Vert U \Vert_F^4.
    \end{aligned}
    \end{equation*}
    Remark that \(U^\top U\) is symmetric and positive semi-definite with eigenvalues \(\lambda_1, \dots, \lambda_n \geq 0\).
    Therefore, \(\trace(U^\top U) = \sum_{i = 1}^n \lambda_i = \Vert U \Vert_F^2\), and
    \begin{equation*}
        \trace(U^\top U U^\top U U^\top U) = \trace((U^\top U)^3) = \sum_{i = 1}^3 \lambda_i^3.
    \end{equation*}
    The power mean inequality (\cref{lem:power-mean}) with \(p = 3\) and \(q = 1\) yields
    \[
        \left( \frac{1}{n} \sum_{i = 1}^n \lambda_i^3 \right)^{\nicefrac{1}{3}} \geq \left( \frac{1}{n} \sum_{i = 1}^n \lambda_i \right)
    \]
    and hence
    \begin{equation*}
        \trace(U^\top U U^\top U U^\top U) = \trace((U^\top U)^3) = \sum_{i = 1}^n \lambda_i^3 \geq \frac{1}{n^2} \left(\sum_{i = 1}^n \lambda_i \right)^3 = \frac{1}{n^2} \Vert U \Vert_F^6.
    \end{equation*}
\end{proof}

\begin{proof}
    We trace the steps from \cref{th:phase-retrieval-smoothness}.
    Since \(f\) is a fourth-order polynomial in the variables \(U\), the Hessian norm can be upper bounded by a second-order polynomial in \(\Vert U \Vert\).
    By combining \cref{lem:matrix-factorization-symmetric-lower-bound,lem:squared-norm-to-norm-lower-bound}, we conclude that the gradient norm is lower bounded by a polynomial of degree \(3\) in \(\Vert U \Vert\), with strictly positive leading coefficient.
    Therefore, \cref{assump:novel-condition} holds.
    The claim then follows from \cref{thm:l0-l1-smoothness-polynomials,thm:anisotropic-smoothness-polynomials}, respectively.
\end{proof}

\subsection{Proof of theorem \ref{th:matrix-factorization-asymmetric-smoothness}}

\begin{lemma}[Gradient norm lower bound] \label{lem:matrix-factorization-asymmetric-lower-bound}
    Consider the asymmetric matrix factorization problem with objective \eqref{eq:asymmetric-factorization-objective}.
    Denote \( V = \max(\Vert W \Vert, \Vert H \Vert_F)\).
    Then,
    \begin{equation*}
        \Vert \nabla f(W, H) \Vert^2 \geq \kappa^2 V^6 - 4 (1+\kappa) \Vert Y \Vert_F V^4.
    \end{equation*}
\end{lemma}
\begin{proof}
    We have that
    \begin{equation*}
    \begin{aligned}
        \Vert \nabla f(W, H) \Vert^2 &= \Vert \nabla_W f(W, H) \Vert_F^2 + \Vert \nabla_H f(W, H) \Vert_F^2.
    \end{aligned}
    \end{equation*}
    The first term can be lower bounded by
    \begin{equation*}
        \begin{aligned}
            \Vert \nabla_W f(W, H) \Vert_F^2 &= \trace((H (H^\top W - Y^\top) + \kappa \Vert W \Vert_F^2 W^\top)((W H - Y) H^\top + \kappa \Vert W \Vert_F^2 W))\\
            &= \trace(H(H^\top W^\top - Y^\top)(W H - Y) H^\top)\\
            &\qquad+ 2 \kappa \Vert W \Vert_F^2 \trace(W^\top (W H - Y) H^\top) + \kappa^2 \Vert W \Vert_F^6\\
            &= \Vert W H H^\top \Vert_F^2 - 2 \trace(H Y^\top W H H^\top) + \trace(H Y^\top Y H^\top)\\
            &\qquad+ 2 \kappa \Vert W \Vert_F^2 \trace(W^\top (W H - Y) H^\top) + \kappa^2 \Vert W \Vert_F^6\\
            &\geq - 2 \Vert Y \Vert_F \Vert W \Vert_F \Vert H \Vert_F^3 + \Vert Y H^\top\Vert_F^2 + 2 \kappa \Vert W \Vert_F^2 (\trace(W^\top W H H^\top)\\
            &\qquad- \trace(W^\top Y H^\top) ) + \kappa^2 \Vert W \Vert_F^6\\
            &\geq - 2 \Vert Y \Vert_F \Vert W \Vert_F \Vert H \Vert_F^3 + 2 \kappa \Vert W \Vert_F^2 \trace(H^\top W^\top W H)\\
            &\qquad- 2 \kappa \Vert W \Vert_F^3 \Vert Y \Vert_F \Vert H \Vert_F + \kappa^2 \Vert W \Vert_F^6\\
            &= - 2 \Vert Y \Vert_F \Vert W \Vert_F \Vert H \Vert_F^3 + 2 \kappa \Vert W \Vert_F^2 \Vert W H \Vert_F^2\\
            &\qquad- 2 \kappa \Vert W \Vert_F^3 \Vert Y \Vert_F \Vert H \Vert_F + \kappa^2 \Vert W \Vert_F^6\\
            &= - 2 \Vert Y \Vert_F \Vert W \Vert_F \Vert H \Vert_F^3 - 2 \kappa \Vert W \Vert_F^3 \Vert Y \Vert_F \Vert H \Vert_F + \kappa^2 \Vert W \Vert_F^6.
        \end{aligned}
    \end{equation*}
    Here the second to last step used the cyclic property of the trace. In a similar way, we obtain
    \begin{equation*}
        \Vert \nabla_H f(W, H) \Vert_F^2 \geq - 2 \Vert Y \Vert_F \Vert H \Vert_F \Vert W \Vert_F^3 - 2 \kappa \Vert H \Vert_F^3 \Vert Y \Vert_F \Vert W \Vert_F + \kappa^2 \Vert H \Vert_F^6.
    \end{equation*}
    Putting this all together, we have that
    \begin{equation*}
        \begin{aligned}
            \Vert \nabla f(W, H) \Vert^2 &\geq \kappa^2 \left( \Vert W \Vert_F^6 + \Vert H \Vert_F^6 \right) - 2 (1+\kappa) \Vert Y \Vert_F \left( \Vert W \Vert_F \Vert H \Vert_F^3 + \Vert W \Vert_F^3 \Vert H \Vert_F \right)\\
            &\geq \kappa^2 V^6 - 4 (1+\kappa) \Vert Y \Vert_F V^4.
        \end{aligned}
    \end{equation*}
\end{proof}

\begin{proof}

    Let \(x\) be the concatenation of \(\vecop(W)\) and \(\vecop(H)\), and let \(V := \max \left(\Vert W \Vert_F, \Vert H \Vert_F\right)\).
    Then \(\Vert x \Vert = \sqrt{\Vert W \Vert_F^2 + \Vert H \Vert_F^2} \geq V \geq \frac{1}{2} \sqrt{\Vert W \Vert_F^2 + \Vert H \Vert_F^2} = \frac{1}{2} \Vert x \Vert\), and we remark that \(\Vert x \Vert \to \infty\) if and only if \(V \to \infty\).

    We now trace the steps from \cref{th:phase-retrieval-smoothness}.
    Since \(f\) is a fourth-order polynomial in the variables \(x\), the Hessian norm can be upper bounded by a second-order polynomial in \(\Vert x \Vert\).
    By combining \cref{lem:matrix-factorization-asymmetric-lower-bound,lem:squared-norm-to-norm-lower-bound} with \(\Vert x \Vert \geq V \geq \frac{1}{2} \Vert x \Vert\), we conclude that the gradient norm is lower bounded by a polynomial of degree \(3\) in \(\Vert x \Vert\), with strictly positive leading coefficient.
    Therefore, \cref{assump:novel-condition} holds.
    The claim then follows from \cref{thm:l0-l1-smoothness-polynomials,thm:anisotropic-smoothness-polynomials}, respectively.
\end{proof}

\subsection{Proof of theorem \ref{th:burer-monteiro-smoothness}}

We first establish the following lemma.
\begin{lemma} \label{lem:burer-monteiro-lower-bound}
    Consider the Burer-Monteiro factorization \eqref{eq:burer-monteiro} of the MaxCut-type SDP \eqref{eq:maxcut} and let \(L_{\beta}\) denote the augmented Lagrangian with penalty parameter \(\beta > 0 \) of this factorized problem.
    Then, there exist constants \(C_1, C_0 \geq 0\) such that
    \begin{equation*}
        \Vert \nabla_x L_{\beta}(x, y) \Vert \geq \frac{2\beta}{n} \Vert x \Vert^3 - C_1 \Vert x \Vert - C_0
    \end{equation*}
\end{lemma}

\begin{proof}
    The gradient of the augmented Lagrangian with respect to \(x\) is
    \begin{equation*}
        \nabla_x L_{\beta}(x, y) = \nabla f(x) + JA^\top (x) y + \nabla_x \left( \frac{\beta}{2} \Vert A(x) \Vert^2 \right)
    \end{equation*}
    and since \(f\) and \(A\) are quadratic in \(x\) we can lower bound its norm by
    \begin{equation*}
        \begin{aligned}
            \Vert \nabla_x L_{\beta}(x, y) \Vert &\geq \beta \bigg \Vert \nabla_x \left( \frac{1}{2} \Vert A(x) \Vert^2 \right) \bigg\Vert - \Vert \nabla f(x) \Vert - \Vert JA(x) \Vert \Vert y \Vert\\
            &\geq \beta \bigg \Vert \nabla_x \left( \frac{1}{2} \Vert A(x) \Vert^2 \right) \bigg\Vert - C_1 \Vert x \Vert - C_0.
        \end{aligned}
    \end{equation*}
    The constraint mapping \(A\) has the particular form \[
        A(x) = \diag(V V^\top) - 1_n = \diag \left( \begin{pmatrix}
            x_1^\top\\
            \vdots\\
            x_n^\top
        \end{pmatrix} \begin{pmatrix}
            x_1 & \dots & x_n
        \end{pmatrix} \right) - 1_n = \begin{pmatrix}
            \Vert x_1 \Vert^2 - 1\\
            \vdots\\
            \Vert x_n \Vert^2 - 1
        \end{pmatrix},
    \]
    and therefore the gradient of the augmenting term equals
    \begin{equation*}
        \nabla_x \left( \frac{1}{2} \Vert A(x) \Vert^2 \right) = 2 \underbrace{\begin{pmatrix}
            \Vert x_1 \Vert^2 x_1\\
            \vdots\\
            \Vert x_n \Vert^2 x_n
        \end{pmatrix}}_{:= y} - 2 x.
    \end{equation*}
    Since \(\Vert y \Vert^2 = \sum_{i = 1}^n \Vert \Vert x_i \Vert^2 x_i \Vert^2 = \sum_{i = 1}^n \Vert x_i \Vert^6\), we can apply the power mean inequality \cref{lem:power-mean} with \(p = 6\) and \(q = 2\).
    This yields
    \begin{equation*}
        \frac{1}{n^{\nicefrac{1}{6}}} \Vert y \Vert^{\nicefrac{2}{6}} = \left( \frac{1}{n} \sum_{i = 1}^n \Vert x_i \Vert^6 \right)^{\nicefrac{1}{6}} \geq \left( \frac{1}{n} \sum_{i = 1}^n \Vert x_i \Vert^{2} \right)^{\nicefrac{1}{2}} = \frac{1}{n^{\nicefrac{1}{2}}} \Vert x \Vert.
    \end{equation*}
    Thus \( \Vert y \Vert^{\nicefrac{1}{3}} \geq \frac{1}{n^{\nicefrac{1}{3}}} \Vert x \Vert \) or \( \Vert y \Vert \geq \frac{1}{n} \Vert x \Vert^3 \).
    Putting everything together, we obtain
    \begin{equation*}
        \begin{aligned}
            \Vert \nabla_x L_{\beta}(x, y) \Vert &\geq \beta \bigg \Vert \nabla_x \left( \frac{1}{2} \Vert A(x) \Vert^2 \right) \bigg\Vert - C_1 \Vert x \Vert - C_0\\
            &\geq 2 \beta \Vert y \Vert - 2 \beta \Vert x \Vert - C_1 \Vert x \Vert - C_0\\
            &\geq \frac{2\beta}{n} \Vert x \Vert^3 - 2 \beta \Vert x \Vert - C_1 \Vert x \Vert - C_0.
        \end{aligned}
    \end{equation*}
    The claim follows by redefining the constant \(C_1\).
\end{proof}

We now present the proof of \cref{th:burer-monteiro-smoothness}.
\begin{proof}
    We again trace the steps from \cref{th:phase-retrieval-smoothness}.
    Since for fixed multipliers \(y\), the augmented Lagrangian \(L_\beta\) is a fourth-order polynomial in the variables \(x\), it follows that the Hessian norm can be upper bounded by a second-order polynomial in \(\Vert x \Vert\).
    Moreover, from \cref{lem:burer-monteiro-lower-bound}, we know that the gradient norm is lower bounded by a polynomial of degree \(3\) in \(\Vert x \Vert\), with strictly positive leading coefficient.
    Therefore, \cref{assump:novel-condition} holds.
    The claim then follows from \cref{thm:l0-l1-smoothness-polynomials,thm:anisotropic-smoothness-polynomials}, respectively.
\end{proof}

\section{Missing proofs of section \ref{sec:saddle-point-avoidance}}
\subsection{Proof of theorem \ref{th:asymptotic-saddle-avoidance-smooth}}

\begin{proof}
    First, under \cref{assumption:phi-star-c2} the iteration map \(
    T_{\gamma, \bar L^{-1}}
    \) is a \(\C^1\) mapping with Jacobian \[
    J T_{\gamma, \bar L^{-1}}(x) = I - \gamma \bar L^{-1} \nabla^2 \phi^*(\bar L^{-1} \nabla f(x)) \nabla^2 f(x).
    \]
    Since \(\lambda_{\max} \left( \nabla^2 \phi^*(\bar L^{-1} \nabla f(x)) \nabla^2 f(x) \right) \leq L \bar L\) and \(\gamma < \frac{1}{L}\) it follows that
    \begin{equation*}
        \lambda_{\min} \left( J T_{\gamma, \bar L^{-1}}(x) \right) \geq 1 - \gamma \bar L^{-1} L \bar L > 0.
    \end{equation*}
    Hence, \( \det J T_{\gamma, \bar L^{-1}}(x) \neq 0\) for all \(x \in \R^n\).
    
    Second, denote the set of unstable fixed points by
    \begin{equation*}
        \cA^\star = \left\{ x^\star \mid T_{\gamma, \bar L^{-1}}(x^\star) = x^\star, \quad \max_i \vert \lambda_i(J T_{\gamma, \bar L^{-1}}(x^\star)) \vert > 1 \right\}.
    \end{equation*}
    For any stationary point \( x^\star\), satisfying \( \nabla f(x^\star) = 0 \), we know that \( x^\star = T_{\gamma, \bar L^{-1}}(x^\star) \), and therefore
    \begin{equation*}
        \begin{aligned}
            \max_i \vert \lambda_i(J T_{\gamma, \bar L^{-1}}(x^\star)) \vert &= \max_i \vert \lambda_i(I - \gamma \bar L^{-1} \nabla^2 \phi^*(\bar L^{-1} \nabla f(x^\star)) \nabla^2 f(x^\star)) \vert\\
            &= \max_i \vert \lambda_i(I - \gamma \bar L^{-1} \nabla^2 f(x^\star)) \vert.
        \end{aligned}
    \end{equation*}
    Hence, if \(\lambda_{\min}(\nabla^2 f(x^\star)) < 0\), then \(\max_i \vert \lambda_i(J T_{\gamma, \bar L^{-1}}(x^\star)) \vert > 1\).
    We conclude that \( \cX^\star \subseteq \cA^\star \).
    
    Thus, we have established all conditions of \cite[Corollary 1]{lee_first-order_2019}, from which the claim follows.
\end{proof}

\subsection{Proof of theorem \ref{th:asymptotic-saddle-avoidance-nonsmooth}}

\begin{proof}
    First, remark that \(U\) is an open set.
    For \(x \in U\) we have by \cite[Lemma 1.3]{oikonomidis_nonlinearly_2025} that
    \begin{equation*}
        \begin{aligned}
            \nabla \phi^*(\bar L^{-1} \nabla f(x)) &= {h^*}'(\Vert \bar L^{-1} \nabla f(x) \Vert) \overline \sign (\bar L^{-1} \nabla f(x))\\
            &= \min \left\{ 1, \Vert \bar L^{-1} \nabla f(x) \Vert \right\} \overline \sign(\bar L^{-1} \nabla f(x))\\
            &= \min \left\{ \frac{1}{\Vert \bar L^{-1} \nabla f(x) \Vert}, 1 \right\} \bar L^{-1} \nabla f(x).
        \end{aligned}
    \end{equation*}
    We partition \( U \) into two sets
    \begin{equation*}
        U_1 := \left\{ x \in \R^n \mid \Vert \bar L^{-1} \nabla f(x) \Vert < 1 \right\}, \quad U_2 := \left\{ x \in \R^n \mid \Vert \bar L^{-1} \nabla f(x) \Vert > 1 \right\}.
    \end{equation*}
    If \( x \in U_1 \), then \(\nabla \phi^*(\bar L^{-1} \nabla f(x)) = \bar L^{-1} \nabla f(x)\).
    Thus, \( T_{\gamma, \bar L^{-1}} \) is continuously differentiable with locally Lipschitz continuous Jacobian on \(U_1\).
    If \( x \in U_2 \), then \(\nabla \phi^*(\bar L^{-1} \nabla f(x)) = \nicefrac{\nabla f(x)}{\Vert \nabla f(x)\Vert} \). 
    Since the function \( \mapsto \nicefrac{z}{\Vert z \Vert} \) is twice continuously differentiable on \( \R^n \setminus \{0\} \) and since \(f \in \cC^{2+}\), it follows that \( T_{\gamma, \bar L^{-1}} \) is continuously differentiable with locally Lipschitz continuous Jacobian on \(U_2\).
    Since \( U = U_1 \cup U_2 \), we conclude that \( J T_{\gamma, \bar L^{-1}} \) is locally Lipschitz continuous on \(U\).    
    By the same arguments as in the proof of \cref{th:asymptotic-saddle-avoidance-smooth} we conclude that
    \begin{equation*}
        \lambda_{\min} \left( J T_{\gamma, \bar L^{-1}}(x) \right) \geq 1 - \gamma \bar L^{-1} L \bar L > 0,
    \end{equation*}
    and that for any \( x^\star \in \cX^\star \) the Jacobian \(J T_{\gamma, \bar L^{-1}}(x^\star)\) is symmetric and has an eigenvalue of absolute value strictly greater than \(1\).
    The claim now follows immediately by applying \cite[Proposition 2.5]{cheridito_gradient_2024}.
\end{proof}

\subsection{Proof of theorem \ref{th:saddle-point-complexity}}

Henceforth, we assume that \(f\) is \((L, \bar L)\)-anisotropically smooth relative to \(\phi\), without further mention.
We first describe the following descent lemma.
\begin{lemma}[Descent lemma] \label{lem:descent-lemma}
    Suppose that \cref{assump:perturbed} holds.
    Let \(\gamma = \nicefrac{\alpha}{L}\), \(\alpha \in (0, 1)\), and let \(\seq{x^k}\) denote a preconditioned gradient sequence, i.e., \(x^{k+1} = T_{\gamma, \lambda}(x^k)\) for \(k \in \N\).
    Then,
    \begin{equation*}
        f(x^{k+1}) - f(x^k) \leq - \frac{\gamma}{2\lambda} \Vert \nabla \phi^*(\lambda \nabla f(x^k)) \Vert^2.
    \end{equation*}
\end{lemma}
\begin{proof}
    By \cite[\S C.1]{oikonomidis_nonlinearly_2025} we have
    \begin{equation*}
        f(x^{k+1}) - f(x^k) \leq - \frac{\alpha}{\lambda L} \phi(\nabla \phi^*(\lambda \nabla f(x^k))).
    \end{equation*}
    The claim then follows from \( \phi(x) \geq \nicefrac{\Vert x \Vert^2}{2} \).
\end{proof}\noindent
We now establish the two key lemmas for our analysis.
First, we establish that if the objective does not decrease much, then the iterates will remain close to the initial point.
\begin{lemma}[Improve or localize] \label{lem:improve-or-localize}
    Suppose that \cref{assump:perturbed} holds.
    Let \(\seq{x^k}\) denote a preconditioned gradient sequence, i.e., \(x^{k+1} = T_{\gamma, \lambda}(x^k)\) for \(k \in \N\).
    Then, for any \( t \geq \tau > 0\)
    \begin{equation*}
        \Vert x^{\tau} - x^0 \Vert \leq \sqrt{2 \lambda \gamma t \left( f(x^0) - f(x^\tau) \right)}.
    \end{equation*}
\end{lemma}\noindent
\begin{proof}
    By consecutively applying the triangle inequality, Cauchy-Schwarz, and \cref{lem:descent-lemma} we obtain
    \begin{equation*}
        \begin{aligned}
            \Vert x^\tau - x^0 \Vert &\leq \sum_{i = 1}^\tau \Vert x^i - x^{i-1} \Vert \leq \sum_{i = 1}^t \Vert x^i - x^{i-1} \Vert \leq \left[ t \sum_{i = 1}^t \Vert x^i - x^{i-1} \Vert^2 \right]^{\nicefrac{1}{2}}\\
            &\leq \left[ \gamma^2 t \sum_{i = 1}^t \Vert \nabla \phi^*(\lambda \nabla f(x^{i-1})) \Vert^2 \right]^{\nicefrac{1}{2}} \leq \sqrt{2 \lambda \gamma t \left( f(x^0) - f(x^t) \right)}.
        \end{aligned}
    \end{equation*}
\end{proof}

Second, we show that the region in which the iterates of \cref{alg:perturbed-GD} remain stuck for at least \(\lceil \cT \rceil\) iterations (if initialized there) is small.
We do this by showing that there exists a point, not to far away, which does yield sufficient decrease.
\begin{lemma}[Coupling sequence] \label{lem:coupling-sequence}
    Suppose that \cref{assumption:hessian-lipschitz,assump:perturbed} holds.
    Let a point \(\tilde x \in \R^n\) satisfy \(\lambda_{\min}(H_\lambda(\tilde x)) \leq - \sqrt{\rho \epsilon}\) where \(\epsilon \leq \frac{L^2}{\rho}\).
    Moreover, let \(\seq{x^k}, \seq{y^k}\) denote two preconditioned gradient sequences, i.e., \(x^{k+1} = T_{\gamma, \lambda}(x^k)\) and \(y^{k+1} = T_{\gamma, \lambda}(y^k)\) for \(k \in \N\), which additionally satisfy 
    \begin{equation*}
        \max \left\{ \Vert x^0 - \tilde x \Vert, \Vert y^0 - \tilde x \Vert \right\} \leq \gamma r, \quad \text{and} \quad x^0 - y^0 = \gamma r_0 e_1,
    \end{equation*}
    where \(e_1\) is the minimum eigenvector of \(H_\lambda(\tilde x)\) and \(r_0 > \omega := 2^{2-\chi} L \cZ\).
    Then,
    \begin{equation*}
        \min \left\{ f(x^{\cT}) - f(x^0), f(y^{\cT}) - f(y^0) \right\} \leq - \cF.
    \end{equation*}
\end{lemma}\noindent
\begin{proof}
    By contradiction, assume that \( \min \left\{ f(x^{\cT}) - f(x^0), f(y^{\cT}) - f(y^0) \right\} > - \cF\).
    \Cref{lem:improve-or-localize} states that for any \(t \leq \cT\)
    \begin{equation} \label{prf:coupling-sequence-contradiction}
        \begin{aligned}
            \max \left\{ \Vert x^t - \tilde x \Vert, \Vert y^t - \tilde x \Vert \right\} &\leq \max \left\{ \Vert x^t - x^0 \Vert, \Vert y^t - y^0 \Vert \right\} + \max \left\{ \Vert x^0 - \tilde x \Vert, \Vert y^0 - \tilde x \Vert \right\}\\
            &\leq \sqrt{2 \lambda \gamma \cT \max \left\{ f(x^0) - f(x^\cT), f(y^0) - f(y^{\cT}) \right\}} + \gamma r\\
            &\leq \sqrt{2 \lambda \gamma \cT \cF} + \gamma r \leq \cZ.
        \end{aligned}
    \end{equation}
    Here the last step follows from \(\epsilon \leq \frac{L^2}{\rho}\).
    Denote by \(z^t := x^t - y^t\) the difference between the two sequences. Then, it follows by the mean value theorem that
    \begin{equation*}
        \begin{aligned}
            z^{t+1} &= z^t - \gamma \left[ \nabla \phi^*(\lambda x^t) - \nabla \phi^*(\lambda y^t) \right] = \left( I - \gamma \cH \right) z^t - \gamma \Delta_t z^t\\
            &= \underbrace{\left( I - \gamma \cH \right)^{t+1} z^0}_{p(t+1)} - \underbrace{\gamma \sum_{i = 0}^{t} \left( I - \gamma \cH \right)^{t-i} \Delta_i z^i}_{q(t+1)},
        \end{aligned}
    \end{equation*}
    where \(\cH := H_\lambda(\tilde x)\) and \(\Delta_t := \int_{0}^{1} \left[ H_\lambda(y^t + \theta(x^t - y^t)) - \cH \right] d\theta \).
    We show by induction that 
    \begin{equation} \label{eq:coupling-sequence-induction}
        \Vert q(t) \Vert \leq \frac{1}{2} \Vert p(t) \Vert, \quad \forall t \in [0, \cT].
    \end{equation}
    For the base case \(t = 0\), the claim holds trivially, since \(
        \Vert q(0) \Vert = 0 \leq \frac{1}{2} \Vert p(0) \Vert.
    \)
    For the induction step, we assume that the claim holds for \(t\) and show that it also holds for \(t + 1\).
    Since \( z^0 \) lies along the minimum eigenvector of \(H_\lambda(\tilde x)\), we have for any \( \tau \leq t \)
    \begin{equation*}
        \Vert z^{\tau} \Vert \leq \Vert p(\tau) \Vert + \Vert q(\tau) \Vert \leq 2 \Vert p(\tau) \Vert = 2 \Vert (I - \gamma \cH)^\tau z^0 \Vert = 2 (1+\gamma \Gamma)^\tau \gamma r_0,
    \end{equation*}
    where \(\Gamma := -\lambda_{\min}(H_\lambda(\tilde x))\).
    By Lipschitz-continuity of \(H_\lambda\) (cfr.\,\cref{assumption:hessian-lipschitz}) we have 
    \[
    \Vert \Delta_t \Vert \leq \int_{0}^{1} \rho \Vert y^t - \tilde x + \theta (x^t - y^t) \Vert d\theta \leq \rho \max \left\{ \Vert x^t - \tilde x \Vert, \Vert y^t - \tilde x \Vert \right\} \leq \rho \cZ.
    \]
    Combined with \(2 \gamma \rho \cZ \cT = \nicefrac{1}{2}\), we obtain
    \begin{equation*}
        \begin{aligned}
            \Vert q(t+1) \Vert &\leq \bigg \Vert \gamma \sum_{i = 0}^{t} \left( I - \gamma \cH \right)^{t - i} \Delta_i z^i \bigg \Vert \leq \gamma \rho \cZ \sum_{i = 0}^{t} \big \Vert \left( I - \gamma \cH \right)^{t - i} \big \Vert \big \Vert z^i \big \Vert\\
            &\leq 2 \gamma \rho \cZ \sum_{i = 0}^{t} (1 + \gamma \Gamma)^{t} \gamma r_0 \leq 2 \gamma \rho \cZ \cT (1+ \gamma \Gamma)^{t} \gamma r_0\\
            &= 2 \gamma \rho \cZ \cT \Vert p(t+1) \Vert \leq \frac{1}{2} \Vert p(t+1) \Vert.
        \end{aligned}
    \end{equation*}
    This completes the proof of \eqref{eq:coupling-sequence-induction}.
    In turn, we conclude that
    \begin{equation*}
        \begin{aligned}
            \max \left\{\Vert x^\cT - \tilde x \Vert, \Vert y^{\cT} - \tilde x \Vert \right\} &\geq \frac{1}{2} \Vert z^\cT \Vert \geq \frac{1}{2} \Vert p(\cT) \Vert - \frac{1}{2} \Vert q(\cT) \Vert \geq \frac{1}{2} \Vert p(\cT) \Vert - \frac{1}{4} \Vert p(\cT) \Vert\\
            &= \frac{1}{4} \Vert p(\cT) \Vert = \frac{(1 + \gamma \Gamma)^\cT \gamma r_0}{4} \geq 2^{\chi - 2} \gamma r_0 > \cZ.
        \end{aligned}
    \end{equation*}
    Here we used \((1+x)^{\nicefrac{1}{x}} \geq 2\) for \(x \in (0, 1]\).
    This contradicts \eqref{prf:coupling-sequence-contradiction} and completes the proof.
\end{proof}

By combining \cref{lem:improve-or-localize} and \cref{lem:coupling-sequence}, we can show that the iterates of \cref{alg:perturbed-GD} will escape from a strict saddle point with high probability.
\begin{lemma}[Escaping strict saddle points] \label{lem:escaping-strict-saddle}
    Suppose that a point \(\tilde x\) satisfies
    \begin{equation*}
        \lambda^{-1} \phi(\nabla \phi^*(\lambda \nabla f(\tilde x))) \leq \frac{\cG^2}{2}, \quad \text{and} \quad \lambda_{\min}(\nabla^2 \phi^*(\lambda \nabla f(\tilde x)) \nabla^2 f(\tilde x)) \leq - \sqrt{\rho \epsilon}
    \end{equation*}
    for \(\epsilon > 0\) small enough.
    Let \(x^0 := \tilde x + \gamma \xi\) where \(\xi\) is sampled uniformly from a ball with radius \(r\), and that \(
    x^{k+1} = x^k - \gamma \nabla \phi^*(\lambda \nabla f(x^k))
    \) for \(k \in \N\).
    Then,
    \begin{equation*}
        \bP \left( f(x^{\cT}) - f(\tilde x) \leq - \frac{\cF}{2} \right) \geq 1 - \frac{L \sqrt{n}}{\sqrt{\rho \epsilon}} \chi^2 2^{8-\chi}.
    \end{equation*}
\end{lemma}\noindent

\begin{proof}
    We define \(\cX_{\text{stuck}}\) as in \cite[Lemma 20]{jin_nonconvex_2021}, i.e.,
    \begin{equation*}
        \cX_{\text{stuck}} := \left\{ x^0 \mid \Vert x^0 - \tilde x \Vert \leq \gamma r \text{ and } f(x^\cT) - f(x^0) > -\cF \text{ where } x^{k+1} = T_{\gamma, \lambda}(x^k) \text{ for } k \in \N \right\}.
    \end{equation*}
    By the same arguments as in \cite[Lemma 20]{jin_nonconvex_2021} we conclude that
    \begin{equation*}
        \bP \left( x^0 \in \cX_{\text{stuck}} \right) \leq \frac{L \sqrt{n}}{\sqrt{\rho \epsilon}} \chi^2 2^{8-\chi}.
    \end{equation*}
    We now proceed by showing that \( f(x^{\cT}) - f(\tilde x) \leq -\frac{\cF}{2} \) if \(x^0 \notin \cX_{\text{stuck}}\).
    By anisotropic smoothness, we have
    \begin{equation*}
        f(x^0) - f(\tilde x) \leq \frac{\gamma}{\lambda} \phi(\xi + \nabla \phi^*(\lambda \nabla f(\tilde x))) = \frac{\gamma}{\lambda} h( \Vert \xi + \nabla \phi^*(\lambda \nabla f(\tilde x))\Vert ).
    \end{equation*}
    By monotonicity of \(h\), the triangle inequality and \cref{lem:upper-bound-preconditioned-gradient-norm}, we have
    \begin{equation*}
        f(x^0) - f(\tilde x) \leq \frac{\gamma}{\lambda} h( \Vert \xi \Vert + \Vert \nabla \phi^*(\lambda \nabla f(\tilde x))\Vert ) \leq \frac{\gamma}{\lambda} h(r + \sqrt{\lambda} \cG) \leq \frac{\gamma}{\lambda} h(2 r).
    \end{equation*}
    Here, the last step uses the bound \( \sqrt{\lambda} \cG \leq r\), which follows from \( \cG = \min \left\{ 1, \frac{1}{\sqrt{\lambda}} \right\} r\).
    For \(r\) sufficiently small -- which holds for \(\epsilon\) sufficiently small -- we can further bound \( h(2r) \leq \frac{5}{8} (2r)^2 = \frac{5}{2} r^2\).
    This yields, again using \( \epsilon \leq \frac{L^2}{\rho} \), that
    \(
        f(x^0) - f(\tilde x) \leq \frac{5}{2} \frac{\gamma}{\lambda} r^2 = \nicefrac{\cF}{2}.
    \)
    We conclude that if \(x^0 \notin \cX_{\text{stuck}}\), then
    \begin{equation*}
        f(x^{\cT}) - f(\tilde x) = \left[ f(x^{\cT}) - f(x^0) \right] + \left[ f(x^0) - f(\tilde x) \right] \leq -\cF + \frac{\cF}{2} = -\frac{\cF}{2}.
    \end{equation*}
\end{proof}

\Cref{th:saddle-point-complexity} now follows by combining \cref{lem:descent-lemma,lem:escaping-strict-saddle}.

\begin{proof}
    Let the total number of iterations of \cref{alg:perturbed-GD} be \[
    \]
    \begin{equation*}
        \begin{aligned}
            T = 8 \max \left\{
                \frac{(f(x^0) - f^*) \cT}{\cF}, \lambda \frac{(f(x^0) - f^*)}{2 \gamma \cG^2}
            \right\} &= \frac{8 L (f(x^0) - f^*)}{\bar L \epsilon^2} \max \left\{ 50 \chi^4, 200 \min \left\{ 1, \sqrt{\bar L} \right\}\chi^3 \right\}\\
            &\leq \frac{8 L (f(x^0) - f^*)}{\bar L \epsilon^2} \max \left\{ 50 \chi^4, 200 \chi^3 \right\}.
        \end{aligned}
    \end{equation*}
    Using \( \max \left\{ 50 \chi^4, 200 \chi^3 \right\} \leq 2^8 \cdot 50 \chi^4 \leq 2^{8 + 7 + 13 + \nicefrac{\chi}{4}} \) for all \(\chi > \nicefrac{1}{4}\) and \(  \), we find
    \begin{equation*}
        T \frac{L \sqrt{n}}{\sqrt{\rho \epsilon}} \chi^2 2^{8 - \chi} \leq \frac{L^2 \sqrt{n}}{\sqrt{\rho \epsilon}} \frac{(f(x^0) - f^*)}{\bar L \epsilon^2} 2^{11} 2^{28 + \nicefrac{\chi}{4}} \chi^2 2^{-\chi}
    \end{equation*}
    and since for \(\chi > \frac{1}{4} \) we have \( \chi^2 2^{-\chi} \leq 2^{10-\nicefrac{\chi}{2}} \) it follows that
    \begin{equation*}
        \begin{aligned}
            T \frac{L \sqrt{n}}{\sqrt{\rho \epsilon}} \chi^2 2^{8 - \chi} &\leq \frac{L^2 \sqrt{n}}{\sqrt{\rho \epsilon}} \frac{(f(x^0) - f^*)}{\bar L \epsilon^2} 2^{39} 2^{-\nicefrac{\chi}{4}}.
        \end{aligned}
    \end{equation*}
    Selecting \(c_{\max} = 2^{-39}\), i.e., \(
        \chi \geq 4 \log_2 \left( 2^{39} \frac{L^2 \sqrt{n} \Delta_f}{ \sqrt{\rho} \bar L \epsilon^{\nicefrac{5}{2}} \delta} \right)
    \)
    we obtain 
    \begin{equation*}
        T \frac{L \sqrt{n}}{\sqrt{\rho \epsilon}} \chi^2 2^{8 - \chi} \leq \delta.
    \end{equation*}
    Observe that \( \cG \leq \epsilon\), such that \(\lambda^{-1} \phi(\nabla \phi^*(\lambda \nabla f(x^k))) \leq \frac{\cG^2}{2}\) implies \(\lambda^{-1} \phi(\nabla \phi^*(\lambda \nabla f(x^k))) \leq \epsilon^2\).
    With probability at least \(1 - \delta\), \cref{alg:perturbed-GD} adds a perturbation at most \( \nicefrac{T}{4\cT}\) times to a point, because by \cref{lem:escaping-strict-saddle} we have
    \begin{equation*}
        f(x^{\cT}) \leq f(x^0) - \frac{T}{4\cT} \frac{\cF}{2} \leq f^*.
    \end{equation*}
    Excluding the iterations which follow within \(\cT\) steps after adding a perturbation, we have at most \( \nicefrac{3T}{4} \) iterations left.
    They either satisfy \(\lambda^{-1} \phi(\nabla \phi^*(\lambda \nabla f(x^k))) \geq \frac{\cG^2}{2}\) or are \( \epsilon\)-second-order stationary points.
    Among these, at most \( \nicefrac{T}{4} \) are not second-order stationary points, because by \cref{lem:descent-lemma} we have
    \begin{equation*}
        f(x^T) \leq f(x^0) - T \frac{\gamma}{2\lambda} \frac{\cG^2}{2} < f^*.
    \end{equation*}
    Therefore, we conclude that at least \( \nicefrac{T}{2}\) iterations of \cref{alg:perturbed-GD} are \(\epsilon\)-second-order stationary points.
\end{proof}

\end{document}